\documentclass[14pt]{article}
\usepackage[margin=1in]{geometry}
\usepackage{tikz}
\usepackage[small,bf,hang]{caption}
\usepackage{float}
\usepackage{amsmath,amsthm,amssymb}
\usepackage{graphicx}
\usepackage{subfig}
\usepackage{tikz}
\usepackage{float}
\usepackage{hyperref}
\hypersetup{colorlinks=true,    citecolor = blue,    linkcolor=black,    filecolor=magenta,    urlcolor=cyan}
\usepackage{enumitem}
\usepackage{mathtools}
\usepackage{multirow}
  
\newtheorem{theorem}{Theorem}
\newtheorem{question}{Question}
\newtheorem{corollary}{Corollary}
\newtheorem{proposition}{Proposition}
\newtheorem{remark}{Remark}
\newtheorem{lemma}{Lemma}
\newtheorem{claim}{Claim}

\newtheorem{operation}{Operation}

\title{Frozen colourings in $2K_2$-free graphs}

\author{Manoj Belavadi
\thanks{Department of Mathematics, Wilfrid Laurier University, Waterloo, ON, Canada, N2L 3C5. Email: \texttt{mbelavadi@wlu.ca}. ORCID: 0000-0002-3153-2339. Research supported by the Natural Sciences and Engineering Research Council of Canada (NSERC) grant RGPIN-2016-06517.}
\and Kathie Cameron
\thanks{Department of Mathematics, Wilfrid Laurier University,
Waterloo, ON, Canada, N2L 3C5. Email: \texttt{kcameron@wlu.ca}. 
ORCID: 0000-0002-0112-2494. Research supported by the Natural Sciences and Engineering Research Council of Canada (NSERC) grant RGPIN-2016-06517.}
\and Elias Hildred
\thanks{Department of Mathematics, Wilfrid Laurier University,
Waterloo, ON, Canada, N2L 3C5. Email: \texttt{hild2190@mylaurier.ca}. ORCID: 0009-0001-9923-1217. Research supported by the Natural Sciences and Engineering Research Council of Canada (NSERC) grant RGPIN-2016-06517.}}

\begin{document}

\maketitle

\begin{abstract}
The \emph{reconfiguration graph of the $k$-colourings} of a graph $G$, denoted $\mathcal{R}_k(G)$, is the graph whose vertices are the $k$-colourings of $G$ and two vertices of $\mathcal{R}_k(G)$ are joined by an edge if the colourings of $G$ they correspond to differ in colour on exactly one vertex. A $k$-colouring of a graph $G$ is called \emph{frozen} if it is an isolated vertex in $\mathcal{R}_k(G)$; in other words, for every vertex $v \in V(G)$, $v$ is adjacent to a vertex of every colour different from its colour. 

A clique partition is a partition of the vertices of a graph into cliques. A clique partition is called a $k$-clique-partition if it contains at most $k$ cliques. Clearly, a $k$-colouring of a graph $G$ corresponds precisely to a $k$-clique-partition of its complement, $\overline{G}$.  
A $k$-clique-partition $\mathcal{Q}$ of a graph $H$ is called \emph{frozen} if for every vertex $v \in V(H)$, $v$ has a non-neighbour in each of the cliques of  $\mathcal{Q}$ other than the one containing $v$.

The cycle on four vertices, $C_4$, is sometimes called the \emph{square}; its complement is called $2K_2$.  

We give several infinite classes of $2K_2$-free graphs with frozen colourings. We give an operation which transforms a $k$-chromatic graph with a frozen $(k+1)$-colouring into a $(k+1)$-chromatic graph with a frozen $(k+2)$-colouring. Our operation preserves being $2K_2$-free. It follows that for all $k \ge 4$, there is a $k$-chromatic $2K_2$-free graph with a frozen $(k+1)$-colouring. We prove these results by studying frozen clique partitions in $C_4$-free graphs.

We say a graph $G$ is \emph{recolourable} if $R_{\ell}(G)$ is connected for all $\ell$ greater than the chromatic number of $G$. We prove that every 3-chromatic $2K_2$-free graph is recolourable. 


\end{abstract}

\section{Introduction}
All graphs in this paper are finite and simple. For a simple graph $G$, the \emph{complement} $\overline{G}$ of $G$ is the simple graph with vertex-set $V(G)$ and where $uv$ is an edge of $\overline{G}$ if and only if $uv$ is not an edge of $G$. 
Let $G$ be a finite simple graph with vertex-set $V(G)$ and edge-set $E(G)$. We use $n = |V(G)|$ to denote the number of vertices of $G$ when the context is clear. 
An \emph{independent set} in a graph $G$ is a set of vertices no two of which are joined by an edge; a \emph{clique} is a set of vertices  every pair of which are joined by an edge. For a positive integer $k$, a \emph{$k$-colouring} of $G$ is a 
partition $\mathcal{C}$ of the vertices into at most $k$ independent sets, called  \emph{colour classes}. A \emph{$k$-clique-partition} is a partition $\mathcal{Q}$ of the vertices into at most $k$ cliques. Clearly, $\mathcal{C}$ is a $k$-colouring of $G$ if and only if $\mathcal{C}$ is a $k$-clique-partition of $\overline{G}$.

We say that $G$ is \emph{$k$-colourable} if it admits a $k$-colouring and is \emph{$q$-clique-partitionable} if it admits a $q$-clique-partition. The \emph{chromatic number} of $G$, denoted $\chi(G)$, is the smallest integer $k$ such that $G$ is $k$-colourable and the \emph{clique partition number} of $G$, denoted $\theta(G)$, is the smallest integer $q$ such that $G$ is $q$-clique-partitionable. Clearly, $\chi(G) = \theta(\overline{G})$. A graph $G$ whose chromatic number is $k$ is called $k$-chromatic.



The \emph{reconfiguration graph of the $k$-colourings}, denoted $\mathcal{R}_k(G)$, is the graph whose vertices are the $k$-colourings of $G$ and two vertices are joined by an edge in $\mathcal{R}_k(G)$ if the colourings they correspond to differ in colour on exactly one vertex. Equivalently, 
two $k$-colourings are adjacent in $\mathcal{R}_k(G)$ if some vertex $v$ can be moved from the part of the partition it is in (that is, from the colour class it is in) to another part, say $U$, of the partition so that the new partition is a colouring. This can be done exactly when $v$ is not adjacent to any vertex of $U$.
We say that $G$ is \emph{k-mixing} if $\mathcal{R}_k(G)$ is connected, and that $G$ is \emph{recolourable} if $G$ is $k$-mixing for all $k>\chi(G)$.

We can also consider the reconfiguration graph of the $q$-clique-partitions of a graph $G$. The vertices of the reconfiguration graph are the $q$-clique-partitions of $G$ and two vertices are joined by an edge in the reconfiguration graph if some vertex $v$ can be moved from the part of the partition it is in (that is, from the clique it is in) to another part, say $U$, of the partition so that the new partition is a clique partition. This can be done exactly when $v$ is adjacent to every vertex of $U$.

Considering colourings and clique partitions as partitions of the vertex-set of a graph, the reconfiguration graph of the $k$-clique-partitions of $\overline{G}$ is precisely $\mathcal{R}_k(G)$. (We comment that normally in mathematics, a partition is thought of as a set of non-empty sets. In reconfiguration of graph colourings, two colourings of a graph are considered different if some vertex has a different colour in the two colourings. So the sets in the partition are really ordered: interchanging the colours of the vertices in two colour classes gives a different colouring. The same concept of order applies to reconfiguration of clique partitions. Also, some of the sets of a colouring or a clique partition can be empty.)

A $k$-colouring of a graph $G$ is called \emph{frozen} if it is an isolated vertex in $\mathcal{R}_k(G)$; in other words, for every vertex $v \in V(G)$, each of the $k$ colours appears in the closed neighbourhood of $v$, or equivalently, if $v$ has a neighbour in each of the colour classes different from the colour class it is in. One way to show that a graph $G$ is not $k$-mixing is to exhibit a frozen $k$-colouring of $G$. Since every $k$-colouring of $K_k$ is frozen, it is common to study $\mathcal{R}_{k+1}(G)$ for a $k$-colourable graph $G$.

A $q$-clique-partition of a graph $G$ is called \emph{frozen} if for every vertex $v \in V(G)$, $v$ has a non-neighbour in each of cliques of the partition different from the  clique it is in. Note that when considering colourings and clique partitions as partitions of the same set $V$ of vertices, a partition corresponding to a colouring of $G$ is frozen if and only if the same partition, considered as a clique partition of $\overline{G}$, is frozen.

Dunbar et al. \cite{dunbar2000} used the term \emph{fall colouring} for frozen colouring, and proved that for each $k \ge 3$, the problem of deciding whether an input graph admits a frozen $k$-colouring is NP-complete.

The cycle on six vertices, $C_6$, admits a frozen 3-colouring, and has the smallest number of vertices of a graph $G$ which admits a frozen $k$-colouring where $k > \chi(G)$. In fact, a cycle $C_n$ admits a frozen 3-colouring if and only if  $n \equiv 0(mod~3)$.

\section{Preliminaries}
\label{sec:pre}

For a vertex $v \in V(G)$, the \emph{open neighbourhood}, $N(v)$, of $v$ is the set of vertices adjacent to $v$ in $G$. The \emph{closed neighbourhood, $N[v]$, of $v$} is the set of vertices adjacent to $v$ in $G$ together with $v$. 

As usual, let $P_n$, $C_n$, and $K_n$ denote the path, cycle, and complete graph on $n$ vertices, respectively. We sometimes refer to $K_3$ as a \emph{triangle} and $C_4$ as a \emph{square}. 

For two vertex-disjoint graphs $G$ and $H$, the \emph{disjoint union} of $G$ 
and $H$, denoted by $G + H$, is the graph with vertex-set $V(G) \cup V(H)$ and edge-set $E(G) \cup E(H)$. For a
positive integer $t$, we use $tG$ to denote the graph obtained from the disjoint union of $t$ copies of G. In particular, the graph $2K_2$ consists of the disjoint union of two copies of $K_2$. The complement of $2K_2$ is $C_4$.  The \emph{paw} is the graph on four vertices consisting of a $K_3$ together with another vertex adjacent to exactly one vertex of the $K_3$. The \emph{diamond} is $K_4$ with one edge deleted (often referred to as $K_4-e$). The edge of the diamond whose end-vertices are of degree 3 is called the \emph{middle edge}.

The subgraph of a graph $G$ \emph{induced} by a subset $S \subseteq  V(G)$ is the graph whose
vertex-set is $S$ and whose edge-set is all edges of $G$ with both ends in $S$.
For a fixed graph $H$, graph $G$ is $H$-free if no induced subgraph of $G$ is isomorphic to $H$. For a set $\mathcal{H}$ of graphs, $G$ is $\mathcal{H}$-free if $G$ is $H$-free for every $H \in \mathcal{H}$. 

A \emph{universal vertex} in a graph $G$ is a vertex which is adjacent to every other vertex of $G$. An \emph{isolated vertex} in a graph $G$ is a vertex which is not adjacent to any vertex of $G$. The \emph{join} of two graphs $G$ and $H$ is obtained by adding all edges between a vertex of $G$ and a vertex of $H$. Two sets of vertices are called \emph{anticomplete (to eachother)} if there is no edge with one end in one set and the other end in the other set. Two sets of vertices are called \emph{complete (to each other)} if there are all possible edges with one end in one set and the other end in the other set.

It is quite easy to see and is used in several papers (see, for example, \cite{feghali2021}) that:

\begin{proposition}
\label{prop:join}
 If $G$ is a $k$-chromatic graph which admits a frozen $\ell$-colouring and if $H$ is an $r$-chromatic graph which admits a frozen $s$-colouring, then the join of $G$ and $H$ is a $(k+r)$-chromatic graph which admits a frozen $(\ell+s)$-colouring. 
\end{proposition}

A \emph{perfect matching} $M$ in a graph $G$ is a set of edges such that each vertex of $G$ is incident to exactly one edge of $M$. For an integer $t \ge 2$, let $K_{t,t}$ denote the complete bipartite graph with $t$ vertices in each part, and let $B_t$ denote $K_{t,t}$ with a perfect matching removed. In \cite{dunbar2000} and \cite{cereceda2008}, it was proved that $B_t$ has a frozen $t$-colouring. Note that $B_t$ is $P_6$-free.

\section{Our contributions}
\label{sec:contributions}
 
A question that has received some attention (see for example, \cite{bonamy2018} and \cite{feghali2021}) is: 

\begin{question}
\label{question}
Given positive integers $k$ and $t$, does there exist a $k$-colourable $P_t$-free graph with a frozen $(k+1)$-colouring? 
\end{question}

The graphs $B_t$ show that for all $t \ge 6$ and $k\ge 2$, the answer to the question is yes. Bonamy and Bousquet \cite{bonamy2018} proved that every $P_4$-free graph $G$ is $k$-mixing for all $k > \chi(G)$, thus for $t \le 4$, the answer to the question is no.

Feghali and Merkel \cite{feghali2021} gave a $7$-chromatic $2K_2$-free graph $G$ on 16 vertices which admits a frozen $8$-colouring. For each positive integer $p$, they then obtained a $7p$-chromatic $2K_2$-free graph which admits a frozen $8p$-colouring by taking $p$ copies of their graph and adding all possible edges between the copies (that is, by taking the pairwise join of $p$ copies of the graph). Thus the answer to Question \ref{question} is yes for $t=5$ and $k\equiv 7(mod~8)$. By adding universal vertices, the result holds for $t=5$ and $k \ge 7$. Feghali and Merkel \cite{feghali2021} asked about the remaining cases. We answer this in the negative for $k \in \{4,5,6\}$ by giving, for all $k \ge 4$, a $k$-chromatic $2K_2$-free graph which admits a frozen $(k+1)$-colouring. Our graphs have the property that their complements are connected (and thus the graphs cannot be decomposed by the join operation). 

We say a graph $G$ is \emph{recolourable} if $R_{\ell}(G)$ is connected for all $\ell\geq \chi(G)$+1. In Section \ref{sec:chi=3}, we prove that every 3-chromatic $2K_2$-free graph is recolourable. Thus the only remaining case of Question \ref{question} is when $t=5$, $k=3$, and the graph contains a $2K_2$.

In \cite{BCM}, it was proved that for a fixed graph $H$, every $H$-free graph is recolourable if and only if $H$ is an induced subgraph of $P_4$ or of $K_3+K_1$. Where $H_1$ and $H_2$ are two fixed graphs on four vertices, it was determined in \cite{BC2024} whether or not all $(H_1,H_2)$-free graphs were recolourable except for $(2K_2, K_4)$-free graphs. This class of graphs is known to be 4-colourable \cite{2K2K4}. Further in \cite{BC2024}, it was proved that every $(2K_2, K_3)$-free graph is recolourable. Thus the result of Section \ref{sec:chi=3} comes close to a dichotomy theorem for recolourability when two graphs on four vertices are forbidden as induced subgraphs. The only open case remaining is whether all 4-chromatic $(2K_2, K_4)$-free graphs which contain a triangle are recolourable.


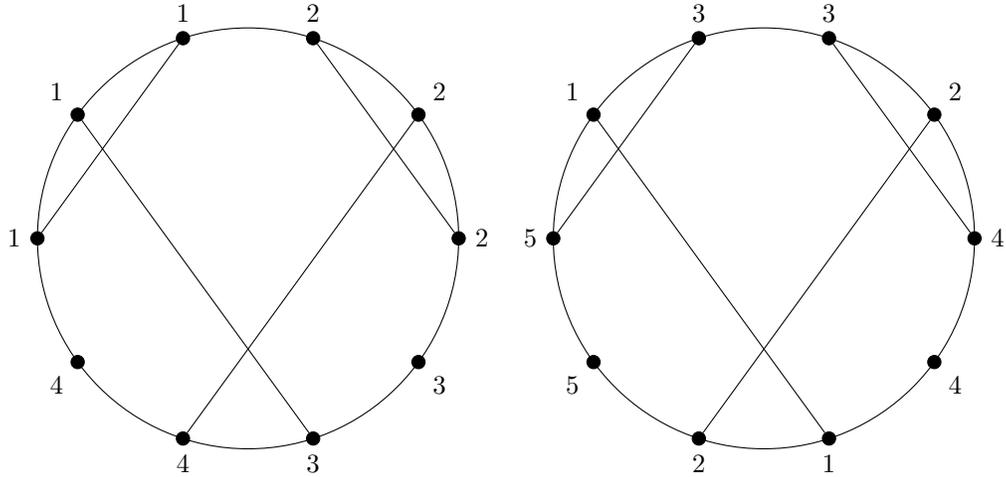
\begin{figure}
\centering
\begin{tikzpicture}[scale=2.8]
\tikzstyle{vertex}=[circle, draw, fill=black, inner sep=0pt, minimum size=5pt]    
        \draw (0,0) circle (1); 
        \node[vertex, label=right:2](0) at (1,0) {};
    \node[vertex, label=above right:2](1) at (cos{36},sin{36}) {};
    \node[vertex, label=above:2](2) at (cos{72},sin{72}) {};
    \node[vertex, label=above:1](3) at (cos{108}, sin{108}) {};
    \node[vertex, label=above left:1](4) at (cos{144}, sin{144}) {};
    \node[vertex, label=left:1](5) at (cos{180}, sin{180}) {};
    \node[vertex, label=below left:4](6) at (cos{216}, sin{216}) {};
    \node[vertex, label=below:4](7) at (cos{252},sin{252}) {};
    \node[vertex, label=below:3](8) at (cos{288}, sin{288}) {};
    \node[vertex, label=below right:3](9) at (cos{324},sin{324}) {};    \draw(0)--(2);\draw(1)--(7);\draw(3)--(5);\draw(4)--(8);
\end{tikzpicture}
\hspace{0mm}
\begin{tikzpicture}[scale=2.8]
\tikzstyle{vertex}=[circle, draw, fill=black, inner sep=0pt, minimum size=5pt]    
        \draw (0,0) circle (1); 
        \node[vertex, label=right:4](0) at (1,0) {};
    \node[vertex, label=above right:2](1) at (cos{36},sin{36}) {};
    \node[vertex, label=above:3](2) at (cos{72},sin{72}) {};
    \node[vertex, label=above:3](3) at (cos{108}, sin{108}) {};
    \node[vertex, label=above left:1](4) at (cos{144}, sin{144}) {};
    \node[vertex, label=left:5](5) at (cos{180}, sin{180}) {};
    \node[vertex, label=below left:5](6) at (cos{216}, sin{216}) {};
    \node[vertex, label=below:2](7) at (cos{252},sin{252}) {};
    \node[vertex, label=below:1](8) at (cos{288}, sin{288}) {};
    \node[vertex, label=below right:4](9) at (cos{324},sin{324}) {};    
    \draw(0)--(2);\draw(1)--(7);\draw(3)--(5);\draw(4)--(8);
     
\end{tikzpicture}

\caption{A square-free graph $\overline{ME_2}$ with a 4-clique-partition (left) and a frozen 5-clique-partition (right). The numbers indicate which clique a vertex is in. Equivalently, the numbers indicate a 4-colouring of the complement $ME_2$ of the graph shown (left) and a frozen 5-colouring of $ME_2$ (right).}
\label{fig:ME2}
\end{figure}

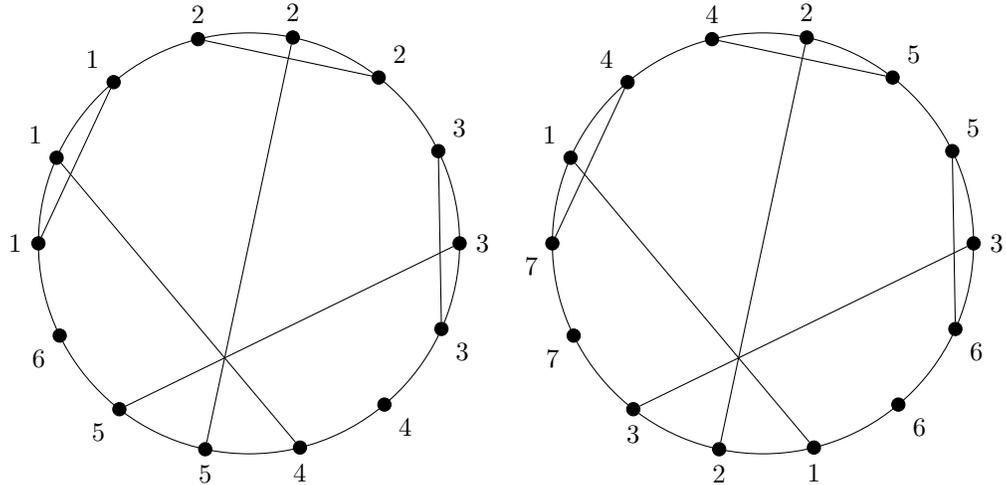
\begin{figure}
\centering
\begin{tikzpicture}[scale=2.8]
\tikzstyle{vertex}=[circle, draw, fill=black, inner sep=0pt, minimum size=5pt]    
        \draw (0,0) circle (1); 
        \node[vertex, label=right:3](0) at (1,0) {};
    \node[vertex, label=above right:3](1) at (cos{26},sin{26}) {};
    \node[vertex, label=above right:2](2) at (cos{52},sin{52}) {};
    \node[vertex, label=above:2](3) at (cos{78},sin{78}) {};
    \node[vertex, label=above:2](4) at (cos{104}, sin{104}) {};
    \node[vertex, label=above left:1](5) at (cos{130}, sin{130}) {};
    \node[vertex, label=above left:1](6) at (cos{156}, sin{156}) {};
    \node[vertex, label=left:1](7) at (cos{180}, sin{180}) {};
    \node[vertex, label=below left:6](8) at (cos{206}, sin{206}) {};
    \node[vertex, label=below left:5](9) at (cos{232},sin{232}) {};
    \node[vertex, label=below:5](10) at (cos{258}, sin{258}) {};
    \node[vertex, label=below:4] (11) at (cos{284},sin{284}) {};    
    \node[vertex, label=below right:4] (12) at (cos{310},sin{310}) {}; 
    \node[vertex, label=below right:3] (13) at (cos{336},sin{336}) {};  
    \draw(0)--(9);\draw(1)--(13);\draw(2)--(4);\draw(3)--(10);
    \draw(5)--(7);\draw(6)--(11);
     
\end{tikzpicture}
\hspace{0mm}
\begin{tikzpicture}[scale=2.8]
\tikzstyle{vertex}=[circle, draw, fill=black, inner sep=0pt, minimum size=5pt]    
        \draw (0,0) circle (1); 
        \node[vertex, label=right:3](0) at (1,0) {};
    \node[vertex, label=above right:5](1) at (cos{26},sin{26}) {};
    \node[vertex, label=above right:5](2) at (cos{52},sin{52}) {};
    \node[vertex, label=above:2](3) at (cos{78},sin{78}) {};
    \node[vertex, label=above:4](4) at (cos{104}, sin{104}) {};
    \node[vertex, label=above left:4](5) at (cos{130}, sin{130}) {};
    \node[vertex, label=above left:1](6) at (cos{156}, sin{156}) {};
    \node[vertex, label=below left:7](7) at (cos{180}, sin{180}) {};
    \node[vertex, label=below left:7](8) at (cos{206}, sin{206}) {};
    \node[vertex, label=below:3](9) at (cos{232},sin{232}) {};
    \node[vertex, label=below:2](10) at (cos{258}, sin{258}) {};
    \node[vertex, label=below:1] (11) at (cos{284},sin{284}) {};    
    \node[vertex, label=below right:6] (12) at (cos{310},sin{310}) {}; 
    \node[vertex, label=below right:6] (13) at (cos{336},sin{336}) {};  
    \draw(0)--(9);\draw(1)--(13);\draw(2)--(4);\draw(3)--(10);
    \draw(5)--(7);\draw(6)--(11);
     
\end{tikzpicture}

\caption{A square-free graph $\overline{ME_3}$ with a 6-clique-partition (left) and a frozen 7-clique-partition (right). The numbers indicate which clique a vertex is in. Equivalently, the numbers indicate a 6-colouring of the complement $ME_3$ of the graph shown (left) and a frozen 7-colouring of $ME_3$ (right).}
\label{fig:ME3}
\end{figure}

The first and third authors did a computer search on all graphs with at most ten vertices to find $k$-colourable $2K_2$-free graphs which admit a frozen $(k+1)$-colouring. Only two graphs were found. One was the graph we call $ME_2$. See Figure \ref{fig:ME2} for the complementary graph, $\overline{ME_2}$. Note that what is shown in Figure \ref{fig:ME2} is actually a 4-clique-partition and a frozen 5-clique-partition of $\overline{ME_2}$; numbers are used to indicate which clique a vertex is in. The other graph they found is one we call $KM_2$, which is $ME_2$ with one edge added (the edge we will later call $u_1u_2$). See Figure \ref{fig:KM2} for $\overline{KM_2}$. 
Both graphs are 4-chromatic and admit a frozen 5-colouring. We believe that these two graphs are the smallest $k$-colourable $2K_2$-free graphs which admit a frozen $(k+1)$-colouring.

In Section \ref{sec:classes}, we give four  infinite classes of $k$-colourable $2K_2$-free graphs which admit frozen $(k+p)$-colourings for various values of $k$ and $p$. These graphs are connected and not decomposable by the join operation. The graphs we construct are dense, so we study clique-partitions in their complements.

In Section \ref{sec:operation}, we give an operation which transforms a $k$-chromatic graph with a frozen $(k+1)$-colouring into a $(k+1)$-chromatic graph with a frozen $(k+2)$-colouring. Note that the operation requires some restrictions on the colouring and frozen colouring. Further, the operation preserves being $2K_2$-free and does not add universal vertices or use the join operation. Again, our approach is to study clique partitions.

In Section \ref{sec:all}, we combine our previously-mentioned  results to show that for all $k \ge 4$, there is a $k$-chromatic $2K_2$-free graph with a frozen $(k+1)$-colouring. This is an improvement on some of the previous examples of frozen colourings since, in these graphs, the gap between the chromatic number of the graph and the number of colours used in a frozen colouring is 1.

In Section \ref{sec:curious}, we make some remarks and mention some open problems.




\section{\texorpdfstring{$2K_2$}\ -free graphs with chromatic number 3 are recolourable}
\label{sec:chi=3}

We use the following result by Bonamy and Bousquet \cite{bonamy2018}.

\begin{lemma}[Renaming Lemma \cite{bonamy2018}]
\label{lem:recolour}
Let $\beta'$ and $\gamma'$ be two $k$-colourings of $G$ that induce the same partition of vertices into colour classes and let $\ell \ge k+1$. Then $\beta'$ can be recoloured into $\gamma'$ in $\mathcal{R}_{\ell}(G)$ by recolouring each vertex at most 2 times.
\end{lemma}

For graph $G$ and a positive integer $k$, we can think of a $k$-colouring of $G$ as a function $\beta  \colon V(G) \to \{1, 2, \ldots, k\}$ such that for each edge $uv \in E(G)$, 
$\beta(u) \neq \beta(v)$. We use $[k]$ to denote ${1,2,\ldots,k}$.

The \emph{diameter} of a graph is the length of a longest shortest path between any two distinct
vertices of the graph. The \emph{$k$-recolouring diameter} of $G$ is the diameter of $\mathcal{R}_k(G)$.

A bipartite graph $G$ is \emph{chordal bipartite} if it does not contain an induced cycle of length more than four. Note that every $2K_2$-free bipartite graph $G$ is a chordal bipartite graph and hence recolourable with $\ell$-recolouring diameter at most 2$n^2$, for all $\ell\ge \chi(G)$+1 \cite{bonamy2014}. This also follows from the fact that every (2$K_2$, triangle)-free graph $G$ is recolourable with $\ell$-recolouring diameter at most 2$n^2$, for all $\ell\ge \chi(G)$+1 \cite{BC2024}. Here we improve the upper bound on the $\ell$-recolouring diameter of 2$K_2$-free bipartite graphs.

\begin{lemma}\label{lem:dominating2K2}
    Let $G$ be a 2$K_2$-free graph. Suppose $V(G)$ can be partitioned into independent sets $A_1, A_2,\dots, A_i$ such that $A_1$ is $($inclusion-wise$)$ maximal. Then for each $j\in \{2,\dots,i\}$, $A_1$ contains a vertex complete to $A_j$.
\end{lemma}
\begin{proof}
    Let $G$ be a 2$K_2$-free graph. Partition $V(G)$ into independent sets $A_1, A_2,\dots, A_i$ such that $A_1$ is (inclusion-wise) maximal. For each $j\in \{2,\dots,i\}$, choose a vertex in $A_1$, say $x_j$, such that $N(x_j)\cap A_j$ is maximized. If $x_j$ is not complete to $A_j$, there is a vertex $y$ in $A_j$ non-adjacent to $x_j$. By the maximality of $A_1$, $y$ has a neighbour $u$ in $A_1$. By the choice of $x_j$, there is a vertex $v$ in $A_j$ adjacent to $x_j$ but non-adjacent to $u$. Then $\{x_j,v,y,u\}$ induces a 2$K_2$, a contradiction. Thus $x_j$ is complete to $A_j$.
\end{proof}

\begin{theorem}\label{thm:2k2bipartiteN}
    Every 2$K_2$-free bipartite graph $G$ is recolourable with $\ell$-recolouring diameter at most 4$n$, for all $\ell\ge \chi(G)$+1.
\end{theorem}
\begin{proof}
    Let $G$ be a 2$K_2$-free bipartite graph. Let $\ell \ge$ 3. Partition $V(G)$ into independent sets $A_1$ and $A_2$ such that $A_1$ is (inclusion-wise) maximal. Given any $\ell$-colouring of $G$ we prove that we can reach a 2-colouring of $G$ that partitions the vertex-set into $A_1$ and $A_2$ by recolouring each vertex at most once. By the Renaming Lemma, there is a path between any two 2-colourings of $G$ that partition the vertex-set into $A_1$ and $A_2$, where each vertex is recoloured at most twice. Starting from any two $\ell$-colourings of $G$, $\beta$ and $\gamma$, we can reach 2-colourings  $\beta^{'}$ and $\gamma^{'}$ in $R_{\ell}(G)$, respectively, which partition the vertex-set into $A_1$ and $A_2$. Then we can obtain $\gamma$ from $\beta$ by recolouring vertices starting from $\beta$ to $\beta^{'}$ to $\gamma^{'}$ to $\gamma$. Each vertex will be recoloured at most 4 times to go from $\beta$ to $\gamma$ in $R_{\ell}(G)$.

    By Lemma \ref{lem:dominating2K2}, $A_1$ contains a vertex, say $x$, complete to $A_2$. Let $\beta$ be any $\ell$-colouring of $G$. There is no vertex in $A_2$ coloured $\beta(x)$. Recolour each vertex in $A_1$ with the colour $\beta(x)$ and recolour each vertex in $A_2$ with a colour $c\neq \beta(x)$. Starting from $\beta$, we have reached a colouring which partitions the vertex-set into $A_1$ and $A_2$, by recolouring each vertex at most once.
\end{proof}

\begin{theorem}
    Every 3-chromatic 2$K_2$-free graph $G$ is recolourable with $\ell$-recolouring diameter at most 14$n$, for all $\ell\ge \chi(G)$+1.
\end{theorem}
\begin{proof}
    Let $G$ be a 3-chromatic 2$K_2$-free graph. Let $\ell\ge$ 4 and let $[\ell]$ be the set of available colours. Partition $V(G)$ into independent sets $A_1$, $A_2$, and $A_3$ such that $A_1$ is (inclusion-wise) maximal. We need some $\chi$-colourings of $G$ to act as anchor points. We say a $\chi$-colouring (i.e., a 3-colouring) of $G$ is \emph{canonical} if it partitions the vertex-set into $A_1$, $A_2$, and $A_3$. By the Renaming Lemma, for all $\ell\ge 4$, there is a path between any two canonical colourings in $R_{\ell}(G)$ where each vertex is recoloured at most twice. Starting from any two $\ell$-colourings of $G$, $\beta$ and $\gamma$, we prove that we can reach canonical colourings  $\beta^{'}$ and $\gamma^{'}$ in $R_{\ell}(G)$, respectively, by recolouring each vertex at most 6 times. Then we can obtain $\gamma$ from $\beta$ by recolouring vertices starting from $\beta$ to $\beta^{'}$ to $\gamma^{'}$ to $\gamma$. Each vertex will be recoloured at most 14 (= 6+2+6) times to go from $\beta$ to $\gamma$ in $R_{\ell}(G)$.

    \begin{claim}\label{C:1partition}
        Any $\ell$-colouring of $G$ which assigns only one colour to some part $A_i$, $i\in [3]$, can be recoloured to a canonical colouring by recolouring each vertex in $V(G)\setminus A_i$ at most 4 times and without recolouring any vertex of $A_i$.
    \end{claim}
    Let $\psi$ be any $\ell$-colouring of $G$ which, for some $i\in [3]$, assigns at most one colour, say $c_i$, to the part $A_i$. Let $A\subseteq V(G)$ be the set of all vertices coloured $c_i$ under $\psi$. Clearly $A_i\subseteq A$. Let $j$ and $k$ be distinct integers in $[3]\setminus \{i\}$. Since $G$-$A$ is a $2K_2$-free bipartite graph, as in the proof of Theorem \ref{thm:2k2bipartiteN}, we can recolour each vertex of $V(G)\setminus A$ at most 4 times to obtain a colouring of $G$ where every vertex of $A_j\setminus A$ is coloured some colour $c_j\neq c_i$ and every vertex of $A_k\setminus A$ is coloured some colour $c_k\notin \{c_i,c_j\}$ without using the colour $c_i$. Recolour each vertex in $A_j\cap A$ with the colour $c_j$ and recolour each vertex in $A_k \cap A$ with the colour $c_k$ to obtain a canonical colouring of $G$. Thus there is a path from $\psi$ to a canonical colouring of $G$ in $R_{\ell}(G)$, for all $\ell\ge \chi(G)$+1, where each vertex of $A_j\cup A_k$ is recoloured at most 4 times.

    \begin{claim}\label{C:dominating}
        If there is a vertex in some $A_i$, $i\in [3]$, adjacent to every vertex not in $A_i$, then any $\ell$-colouring of $G$ can be recoloured to a canonical colouring by recolouring each vertex at most 4 times.
    \end{claim}
    For some $i\in [3]$, let $x$ in $A_i$ be adjacent to every vertex outside $A_i$. Let $\psi$ be any $\ell$-colouring of $G$. Recolour each vertex in $A_i$ with the colour $\psi(x)$. Now, by Claim \ref{C:1partition}, we can reach a canonical colouring of $G$ by recolouring each vertex in $V(G)\setminus A_i$ at most 4 times and without recolouring any vertex in $A_i$. Therefore, we can reach a canonical colouring of $G$ by recolouring each vertex at most 4 times.\vspace{5mm}

    By Lemma \ref{lem:dominating2K2}, there are vertices $x_2$ and $x_3$ in $A_1$ complete to $A_2$ and $A_3$, respectively. By Claim \ref{C:dominating}, we may assume that $x_2$ and $x_3$ are distinct. Let $\beta$ be any $\ell$-colouring of $G$.\vspace{5mm}
    
    Suppose $\beta(x_2)$ = $\beta(x_3)$ = $c_1$, then there is no vertex outside $A_1$ coloured $c_1$. Recolour each vertex in $A_1$ with colour $c_1$. Now, by Claim \ref{C:1partition}, we can reach a canonical colouring of $G$ by recolouring each vertex in $A_2\cup A_3$ at most 4 times and without recolouring any vertex in $A_1$. Therefore, we can reach a canonical colouring of $G$ by recolouring each vertex at most 4 times.\vspace{5mm}

    Suppose $\beta(x_2) \neq \beta(x_3)$. Let $\beta(x_2)$ = 1 and let $\beta(x_3)$ = 2. Note that no vertex of $A_2$ received colour 1 and no vertex of $A_3$ received colour 2. Recolour as many vertices as possible in $A_2$ with colour 2; that is, recolour with colour 2 every vertex of $A_2$ which does not have a neighbour of colour 2 in $A_1$. Recolour as many vertices as possible in $A_3$ with colour 1; that is, recolour with colour 1 every vertex of $A_3$ which does not have a neighbour of colour 1 in $A_1$. Recolour as many vertices as possible in $A_1$ with either colour 1 or 2; that is, for vertex $v$ of $A_1$ which is non-adjacent to a vertex coloured 1 or 2, recolour $v$ with colour 1 if $v$ is non-adjacent to a vertex coloured 1 in $A_3$ or recolour $v$ with colour 2 if $v$ is non-adjacent to a vertex coloured 2 in $A_2$. This new colouring, say $\zeta$, is obtained from $\beta$ by recolouring each vertex at most once.
    
    Now a vertex in $A_1$ is coloured neither colour 1 nor colour 2 if and only if it is adjacent to a vertex coloured 1 in $A_3$ and adjacent to a vertex coloured 2 in $A_2$.
    \begin{claim}\label{C:1colour}
        If there is a vertex in $A_1$ coloured $c\in \{3,4\}$ under $\zeta$, then there are no vertices outside $A_1$ coloured $c$ under $\zeta$.
    \end{claim}
    We prove the claim for $c$ = 3. Let $x\in A_1$ and $y\in A_2\cup A_3$ be coloured 3 under $\zeta$. The vertex $x$ was not recoloured with colour either 1 or 2, because it is adjacent to a vertex $u$ coloured 2 in $A_2$ and adjacent to a vertex $v$ coloured 1 in $A_3$. If $y\in A_2$, then by the choice of $\zeta$, it is adjacent to a vertex $w$ coloured 2 in $A_1$. Then $\{y,w,x,u\}$ induces a 2$K_2$, a contradiction. The proof is similar if $y$ is in $A_3$. This proves Claim \ref{C:1colour}.\vspace{5mm}

    We have two cases.
    
    \noindent
    \textit{Case 1}: There is a vertex in $A_1$ coloured either 3 or 4 under $\zeta$.\\
    Let there be a vertex $x$ coloured either 3 or 4 in $A_1$. Then by Claim \ref{C:1colour} there are no vertices coloured $\zeta(x)$ outside $A_1$. Recolour each vertex in $A_1$ with the colour $\zeta(x)$. Now, by Claim \ref{C:1partition}, we can reach a canonical colouring of $G$ by recolouring each vertex in $A_2\cup A_3$ at most 4 times and without recolouring any vertex in $A_1$. Therefore, starting from $\beta$ we recoloured each vertex at most once to reach $\zeta$ and then recoloured each vertex at most 4 times to reach a canonical colouring of $G$. This completes the proof for Case 1.\vspace{5mm}

    \noindent
    \textit{Case 2}: There is no vertex in $A_1$ coloured either 3 or 4 under $\zeta$.\vspace{5mm}

    \noindent
    \textit{Case 2 $($a$)$}: Suppose for $j\in \{2,3\}$, there is no vertex coloured either 3 or 4 in $A_j$.
    
    Recolour each vertex in $A_i$, $i\in \{2,3\}\setminus \{j\}$, with colour either 3 or 4, respectively. Now, by Claim \ref{C:1partition}, we can reach a canonical colouring of $G$ by recolouring each vertex in $V(G)\setminus A_i$ at most 4 times and without recolouring any vertex in $A_i$. Therefore, starting from $\beta$ we recoloured each vertex at most once to reach $\zeta$ and then recoloured each vertex at most 4 times to reach a canonical colouring of $G$. This completes the proof for Case 2(a).\vspace{5mm}

    \noindent
    \textit{Case 2 $($b$)$}: Both colours 3 and 4 appear on $A_j$ under $\zeta$, for all $j\in \{2,3\}$. 
    
    Suppose there are two vertices $u$ and $v$ in $A_2$ coloured 3 and 4, respectively, such that $u$ has a neighbour $u^{'}$ coloured 4 and $v$ has a neighbour $v^{'}$ coloured 3. Then $u^{'}$ and $v^{'}$ must be in $A_3$. This implies that $\{u, u^{'}, v, v^{'}\}$ induces a 2$K_2$, a contradiction. Therefore there are no two vertices $u$ and $v$ coloured 3 and 4, respectively, in $A_2$ such that $u$ has a neighbour coloured 4 and $v$ has a neighbour coloured 3. Without loss of generality, assume that there is no vertex in $A_2$ coloured 4 which is a adjacent to a vertex coloured 3.

    Recolour each vertex coloured 4 in $A_2$ with colour 3. Now there is no vertex in $A_1\cup A_2$ coloured 4. Recolour each vertex in $A_3$ with colour 4. Now, by Claim \ref{C:1partition}, we can reach a canonical colouring of $G$ by recolouring each vertex in $A_1\cup A_2$ at most 4 times and without recolouring any vertex in $A_3$. Therefore, starting from $\beta$ we recoloured each vertex at most once to reach $\zeta$ and then recoloured each vertex at most 5 times to reach a canonical colouring of $G$. This completes the proof for Case 2(b).
\end{proof}

\section{Four infinite classes of \texorpdfstring{$2K_2$}\ -free graphs which admit frozen colourings}
\label{sec:classes}

A \emph{Hamiltonian cycle} in a graph $G$ is a cycle which contains all the vertices of $G$.\\

For an integer $q \ge 2$, $\overline{{ME}_q}$  is the graph with $4q+2$ vertices \\ 
$\{u_i: i=0,1,\ldots,q+1\} \cup \ \{\cup \{v_{i1}, v_{i2}, v_{i3}\}: i=1,2,\ldots,q\}$\\ whose edges are: 
\begin{itemize}
\item  the edges of a Hamiltonian cycle $C$: $u_0, u_1, \ldots, u_{q+1}, v_{11}, v_{12}, v_{13}, v_{21}, v_{22}, v_{33},\ldots, v_{q1}, v_{q2}, v_{q3}, u_0$
\item edges $u_i v_{i2}$ for $i = 1, 2, \ldots, q$
\item edges $v_{i1}v_{i3}$ for $i = 1, 2, \ldots, q$
\end{itemize}

See Figure \ref{fig:ME2} for $\overline{ME_2}$  and Figure \ref{fig:ME3} for $\overline{ME_3}$.

We refer to  $\{v_{i1}, v_{i2}, v_{i3}\}$ as \emph{triangle $i$}. Note that $\overline{{ME}_q}$ consists of a Hamiltonian cycle $C$ together with $q$ edges which induce $q$ vertex-disjoint triangles with consecutive pairs of edges of $C$, and $q$ more edges $u_i v_{i2}$ each of which induces a paw with triangle $i$. Also note that the only neighbours of vertices $u_0$ and $u_{q+1}$ are their neighbours on $C$. The number of edges of $\overline{{ME}_q}$  is $(4q+2)+2q = 6q+2$. 

\begin{theorem}\label{C4free} 
For $q \ge 2$, $\overline{{ME}_q}$ is $C_4$-free.
\end{theorem}

\begin{proof}
Consider the graph $\overline{{ME}_q}$ where $q \ge 2$. Edge $v_{i1}v_{i3}$ cannot be part of an induced 4-cycle in $\overline{{ME}_q}$ because clearly $v_{i2}$ can't be part of such a cycle and $v_{i1}$'s only other neighbour is either $v_{i-1 \mkern3mu 3}$ if $i \ge 2$ or $u_{q+1}$ if $i=1$, and $v_{i3}$'s  only other neighbour is either $v_{i+1\mkern3mu 1}$ if $i \le q-1$ or $u_0$ if $i=q$, and these neighbours are not adjacent. 

Edge $u_iv_{i2}$ makes two cycles with $C$. The two cycles are generally not induced cycles; a shorter cycle can be obtained by replacing any occurrence of $v_{j1}, v_{j2}, v_{j3}$ by $v_{j1}, v_{j3}$. We first consider cycles containing only one edge of the type $u_iv_{i2}$. The shortest such cycles occur when $i=1$ or $i=q$, and are $v_{12}, u_1, u_2, \ldots, u_{q+1}, v_{11}, v_{12}$  and $v_{q2}, v_{q3}, u_0, u_1, \ldots, u_q, v_{q2}$, respectively, and each has length $q+3 \ge 5$. 

A shortest cycle containing two edges $u_iv_{i2}$ and $u_jv_{j2}$ is when $j=i+1$, and is the 6-cycle: \\$u_i, v_{i2}, v_{i3}, v_{i+1\mkern3mu1}, v_{i+1\mkern3mu2},u_{i+1}, u_i$. 

Thus $\overline{{ME}_q}$ is $C_4$-free.

\end{proof}

\begin{corollary}\label{$2K_2$free}
For $q \ge 2$, ${ME}_q$ is $2K_2$-free.
\end{corollary}


For a graph $G$, $\alpha(G)$ denotes the size of a largest independent set in $G$ and $\omega(G)$ denotes the size of a largest clique in $G$.

\begin{theorem}\label{theta}
\[\text{For } q \ge 2,~ \theta(\overline{{ME}_q})\ =\ \alpha(\overline{{ME}_q})\ =\
\begin{dcases*}
    (3q+2)/2 & \textit{if } q~ \textit{is even} \\
    (3q+3)/2 & \textit{if } q~ \textit{is odd}
\end{dcases*}
\]
\end{theorem}

\begin{proof}
Let $q \ge 2$ be even. Create a clique partition of $\overline{{ME}_q}$ consisting of the following cliques:

\begin{itemize}
\item For $i=1, 2, \dots, q$, let the vertices of triangle $i$ be a clique in the clique partition
\item Divide the vertices of the path from $u_0$ to $u_{q+1}$ into $(q+2)/2$ cliques as follows: $\{u_0,u_1\}, \{u_2,u_3\}, \ldots, \{u_q,u_{q+1}\}$ and put these cliques into the clique partition
\end{itemize}

Vertices $u_0, u_2, \dots, u_q,v_{11}, v_{21}, v_{31}, \ldots, v_{q1}$ form an independent set of size $(q+2)/2 + q$ in $\overline{{ME}_q}$. 

Now let $q \ge 3$ be odd. The proof is similar to the even case, except that the path from $u_0$ to $u_{q+1}$ in $\overline{{ME}_q}$ has an odd number of vertices and thus requires $(q+1)/2 + 1$ colours in ${ME}_q$. Create a clique partition of $\overline{{ME}_q}$ consisting of the following cliques:

\begin{itemize}
    \item For $i=1, 2, \dots, q$, let the vertices of triangle $i$ be a clique in the clique partition 
    \item Divide the vertices of the path from $u_0$ to $u_{q+1}$ into $(q+3)/2$ cliques as follows:\\ $\{u_0,u_1\}, \{u_2,u_3\}, \ldots, \{u_{q-1},u_q\}, \{u_{q+1}\}$ and put these cliques into the clique partition
\end{itemize}

Vertices $u_0, u_2, \dots, u_{q+1}, v_{12}, v_{21}, v_{31}, \ldots, v_{q1}$ form an independent set of size $(q+1)/2 + 1 + q$ in $\overline{{ME}_q}$.

\end{proof}

\begin{corollary}\label{chi}
\[\text{For } q \ge 2, ~ \chi({ME}_q)~ =~ \omega({ME}_q) ~=~
\begin{dcases*}
    (3$q$+2)/2 & \textit{if } q~ \textit{is even} \\
    (3$q$+3)/2 & \textit{if } q~ \textit{is odd}
\end{dcases*}
\]
\end{corollary}

\begin{lemma} \label{lem:triangles}
Let $\mathcal{Q}$ be a partition of the vertex-set $V(G)$ of graph $G$ into cliques of size 2. Then $\mathcal{Q}$ is a frozen clique partition if and only if every triangle of $G$ intersects three distinct cliques of $\mathcal{Q}$.
\end{lemma}

\begin{proof}
Let $\mathcal{Q}$ be a partition of the vertex-set $V(G)$ of graph $G$ into cliques of size 2. Then every triangle of $G$ intersects at least two cliques of $\mathcal{Q}$. 

By definition, $\mathcal{Q}$ is not a frozen clique partition if and only if there is some clique $Q=\{q_1,q_2\} \in \mathcal{Q}$ and some vertex $v \notin Q$ such that $v$ is adjacent to both $q_1$ and $q_2$, which means that triangle $\{v, q_1, q_2\}$ intersects exactly two cliques of $\mathcal{Q}$, namely $Q$ and the clique containing $v$.
\end{proof}

\begin{theorem}\label{thm:frozencp}
For $q \ge 2$, $\overline{{ME}_q}$ has a frozen $(2q+1)$-clique-partition.
\end{theorem}

\begin{proof}
Create a clique partition $\mathcal{Q^*}$ of $\overline{{ME}_q}$ consisting of the following cliques:
\begin{itemize}
\item For $i=1, 2, \dots, q$, let $\{u_i, v_{i2}\}$ be a clique of the clique partition. 
\item For $i=1, 2, \dots, q-1$, let $\{v_{i3}, v_{i+1\mkern3mu1}\}$ be a clique of the clique partition.
\item Let $\{v_{q3},u_0\}$ be a clique of the clique partition.
\item Let $\{u_{q+1},v_{11}\}$ be a clique of the clique partition.
\end{itemize}

In $\overline{{ME}_q}$, the only triangles are triangles 1 to q. It is easily seen that each triangle $i$ intersects three different cliques of $\mathcal{Q^*}$. Thus the result follows from Lemma \ref{lem:triangles}.
\end{proof}

\begin{corollary}\label{cor:frozencol}
For $q \ge 2$, ${ME}_q$ has a frozen $(2q+1)$-colouring.
\end{corollary}

See Table \ref{tab:MEq complements} for parameters of $\overline{{ME}_q}$ graphs and Table \ref{tab:MEq graphs} for parameters of ${{ME}_q}$ graphs. 

\begin{table} 
    \centering
    \begin{tabular}{|c|c|c|c|c|c|c|c|}
    \hline
        $q$ & $n$ & min  & max & \# edges & $\theta=\alpha$ & \# cliques in & (\# cliques in  \\
         &  & degree & degree & &  & frozen clique ptn & frozen clique ptn) - $\theta$ 
         \\
         \hline
         & & & & & &  &     \\
         $q$&$4q+2$ & $2$ & $3$ & $6q+2$&  $(3q+2)/2$ for even $q$  & $n/2=2q+1$ &  $q/2$ for even $q$   \\
         & & & & & $(3q+3)/2$ for odd $q$  &  & $(q-1)/2$ for odd $q$    \\
        & & & &  &  &  &   \\
         \hline
          & & & &  &  & &    \\
         2 & 10 & 2 & 3 & 14  &  4 & 5 & 1  \\
         3 & 14 & 2 & 3 & 20  & 6 & 7 & 1  \\
         4 & 18 & 2 & 3 & 26 & 7 & 9 & 2 \\
         5 & 22 & 2 & 3 & 32 & 9 & 11& 2  \\
         6 & 26 & 2 & 3 & 38 & 10 & 13 & 3 \\
         7 & 30 & 2 & 3 & 44 & 12 & 15 & 3  \\
         8 & 34 & 2 & 3 & 48 & 13 & 17 & 4  \\
         \hline
          
    \end{tabular}
    \caption{Parameters of  $\overline{{ME}_q}$ graphs}
    \label{tab:MEq complements}
\end{table}

\begin{table} 
    \centering
    \begin{tabular}{|c|c|c|c|c|c|c|c|}
    \hline
        $q$ & $n$ & min  & max & \# edges & $\chi=\omega$ & \# colours in & (\# colours in  \\
         &  & degree & degree & &  & frozen colouring & frozen colouring) - $\chi$ 
         \\
         \hline
         & & & & & &  &     \\
         $q$&$4q+2$ & $4q-2$ & $4q-1$ & $8q^2-1$&  $(3q+2)/2$ for even $q$  & $n/2=2q+1$ &  $q/2$ for even $q$   \\
         & & & & & $(3q+3)/2$ for odd $q$  &  & $(q-1)/2$ for odd $q$    \\
        & & & &  &  &  &   \\
         \hline
          & & & &  &  & &    \\
         2 & 10 & ~6 & ~7 & 31  &  4 & 5 & 1  \\
         3 & 14 & 10 & 11 & 71  & 6 & 7 & 1  \\
         4 & 18 & 14 & 15 & 127 & 7 & 9 & 2 \\
         5 & 22 & 18 & 19 & 199 & 9 & 11& 2  \\
         6 & 26 & 22 & 23 & 287 & 10 & 13 & 3 \\
         7 & 30 & 26 & 27 & 391 & 12 & 15 & 3  \\
         8 & 34 & 30 & 31 & 511 & 13 & 17 & 4  \\
         \hline
          
    \end{tabular}
    \caption{Parameters of ${ME}_q$ graphs}
    \label{tab:MEq graphs}
\end{table}

As noted in Section \ref{sec:contributions}, the graph $ME_2$ was found by a computer search as was the graph we will call 
$KM_2$
which is $ME_2$ with edge $u_1u_2$ added. 

We now define a second class of graphs, ${ME}^*_q$ where $q \ge 2$. We obtain ${ME}^*_q$ from ${ME}_q$ by deleting the edge $u_0u_{q+1}$. Equivalently, we obtain $\overline{{ME}^*_q}$ from $\overline{{ME}_q}$ by adding the edge $u_0u_{q+1}$. 

\begin{theorem}\label{C4free-} 
For $q \ge 3$, $\overline{{ME}^*_q}$ is $C_4$-free.
\end{theorem}

\begin{proof}
By Theorem \ref{C4free}, we only have to consider cycles in $\overline{{ME}^*_q}$, $q \ge 3$ containing edge $u_0u_{q+1}$. This edge creates two cycles with $C$, one of which has length $q+2 \ge 5$, and the other which is not induced but a shorter cycle can be obtained from it by using edges $v_{i1}v_{i3}$, so by an argument similar to that given in the proof of Theorem \ref{C4free}, the induced cycle is not a 4-cycle (it has length $2q+2$). There are two 5-cycles  containing $u_0u_{q+1}$: $u_0, u_1, v_{12}, v_{11}, u_{q+1},u_0$ and $u_0, v_{q3}, v_{q2}, u_q, u_{q+1}, u_0$.   
\end{proof}

\begin{remark}
By deleting the edge $u_0u_{q+1} = u_0u_3$ from $ME_2$, we will obtain a $2K_2$ induced by $\{u_0,u_1,u_2, u_3\}.$ The resulting graph $ME^*_2$ is thus not $2K_2$-free, but is $P_5$-free. 
\end{remark}

\begin{corollary}\label{$2K_2$free-}
For $q \ge 3$, ${ME}^*_q$ is $2K_2$-free.
\end{corollary}

\begin{theorem}\label{theta-}
\[\text{For } q \ge 2,~ \theta({ME}^*_q) ~=~ \alpha({ME}^*_q)~ =
  \begin{dcases*}
    (3q+2)/2 & \textit{if } q~ \textit{is even} \\
    (3q+3)/2 & \textit{if } q~ \textit{is odd} \\
  \end{dcases*}
\]
\end{theorem}

\begin{proof}
Graph $\overline{{ME}^*_q}$ is $\overline{{ME}_q}$ with an edge added, so a clique partition of $\overline{{ME}_q}$ is 
a clique partition of $\overline{{ME}^*_q}$. When $q$ is even, the independent set of $\overline{{ME}_q}$ given in the proof of Theorem \ref{theta} is an independent set in $\overline{{ME}^*_q}$ because it does not contain $u_{q+1}$. When $q$ is odd, the clique partition of $\overline{{ME}^*_q}$ can be seen to be minimum because the vertex-set of $\overline{{ME}^*_q}$ can be partitioned into $q$ triangles and one induced odd cycle $C_{q+2}$; the odd cycle $C_{q+2}$ requires at least $(q+3)/2$ cliques in any clique partition. Thus the size of a smallest clique partition of $\overline{{ME}^*_q}$ when q is odd is 
$q+(q+3)/2 = (3q+3)/2$.
\end{proof}

\begin{corollary} \label{chi-}
\[\text{For } q \ge 2,\
\chi({ME}^*_q)\ =\ \omega({ME}^*_q)\ =\
  \begin{dcases*}
    (3$q$+2)/2 & \textit{if } q~ \textit{is even} \\
    (3$q$+3)/2 & \textit{if } q~ \textit{is odd} \\
  \end{dcases*}
\]  
\end{corollary}

\begin{theorem}\label{thm:frozen-}
For $q \ge 2$, $\overline{{ME}^*_q}$ has a frozen $(2q+1)$-clique-partition.
\end{theorem}

\begin{proof}
We claim that the frozen clique partition of $\overline{{ME}_q}$ given in the proof of Theorem \ref{thm:frozencp} is a frozen clique partition of $\overline{{ME}^*_q}$. Note that $\overline{{ME}_q}$ and $\overline{{ME}^*_q}$ have exactly the same set of triangles, so the result follows from Theorem \ref{thm:frozencp} and Lemma \ref{lem:triangles}.
\end{proof}

\begin{corollary}\label{cor:frozen-}
For $q \ge 2$, ${ME}^*_q$ has a frozen $(2q+1)$-colouring.
\end{corollary}


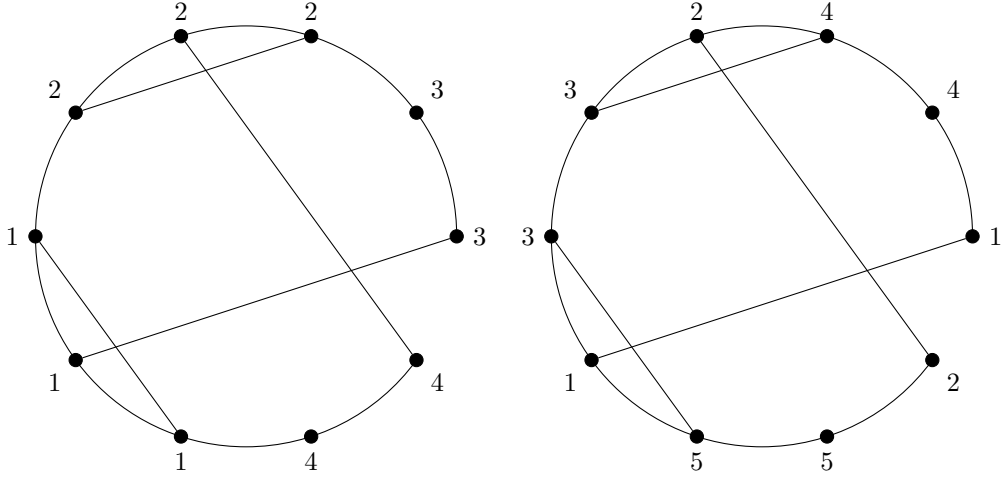
\begin{figure}
\centering
\begin{tikzpicture}[scale=2.8]
\tikzstyle{vertex}=[circle, draw, fill=black, inner sep=0pt, minimum size=5pt]    
        \node[vertex, label=right:3](0) at (1,0) {};
    \node[vertex, label=above right:3](1) at (cos{36},sin{36}) {};
    \node[vertex, label=above:2](2) at (cos{72},sin{72}) {};
    \node[vertex, label=above:2](3) at (cos{108}, sin{108}) {};
    \node[vertex, label=above left:2](4) at (cos{144}, sin{144}) {};
    \node[vertex, label=left:1](5) at (cos{180}, sin{180}) {};
    \node[vertex, label=below left:1](6) at (cos{216}, sin{216}) {};
    \node[vertex, label=below:1](7) at (cos{252},sin{252}) {};
    \node[vertex, label=below:4](8) at (cos{288}, sin{288}) {};
    \node[vertex, label=below right:4](9) at (cos{324},sin{324}) {};    
    \draw (0) arc(0:324:1);
    \draw(0)--(6);\draw(2)--(4);\draw(3)--(9);\draw(5)--(7);
\end{tikzpicture}
\hspace{0mm}
\begin{tikzpicture}[scale=2.8]
\tikzstyle{vertex}=[circle, draw, fill=black, inner sep=0pt, minimum size=5pt]    
        \node[vertex, label=right:1](0) at (1,0) {};
    \node[vertex, label=above right:4](1) at (cos{36},sin{36}) {};
    \node[vertex, label=above:4](2) at (cos{72},sin{72}) {};
    \node[vertex, label=above:2](3) at (cos{108}, sin{108}) {};
    \node[vertex, label=above left:3](4) at (cos{144}, sin{144}) {};
    \node[vertex, label=left:3](5) at (cos{180}, sin{180}) {};
    \node[vertex, label=below left:1](6) at (cos{216}, sin{216}) {};
    \node[vertex, label=below:5](7) at (cos{252},sin{252}) {};
    \node[vertex, label=below:5](8) at (cos{288}, sin{288}) {};
    \node[vertex, label=below right:2](9) at (cos{324},sin{324}) {};    
    \draw (0) arc(0:324:1);
    \draw(0)--(6);\draw(2)--(4);\draw(3)--(9);\draw(5)--(7);
     
\end{tikzpicture}

\caption{A square-free graph $\overline{KM_2}$ with a 4-clique-partition (left) and a frozen 5-clique-partition (right). Equivalently, a 4-colouring of the complement $KM_2$ of the graph shown (left) and a frozen 5-colouring of $KM_2$ (right).}
\label{fig:KM2}
\end{figure}

\begin{figure}
\centering
\begin{tikzpicture}[scale=2.8]
\tikzstyle{vertex}=[circle, draw, fill=black, inner sep=0pt, minimum size=5pt]    
        \node[vertex, label=right:4](0) at (1,0) {};
    \node[vertex, label=above right:4](1) at (cos{26},sin{26}) {};
    \node[vertex, label=above right:3](2) at (cos{52},sin{52}) {};
    \node[vertex, label=above:3](3) at (cos{78},sin{78}) {};
    \node[vertex, label=above:3](4) at (cos{104}, sin{104}) {};
    \node[vertex, label=above left:2](5) at (cos{130}, sin{130}) {};
    \node[vertex, label=above left:2](6) at (cos{156}, sin{156}) {};
    \node[vertex, label=left:2](7) at (cos{180}, sin{180}) {};
    \node[vertex, label=below left:1](8) at (cos{206}, sin{206}) {};
    \node[vertex, label=below left:1](9) at (cos{232},sin{232}) {};
    \node[vertex, label=below:1](10) at (cos{258}, sin{258}) {};
    \node[vertex, label=below:6] (11) at (cos{284},sin{284}) {};    
    \node[vertex, label=below right:6] (12) at (cos{310},sin{310}) {}; 
    \node[vertex, label=below right:5] (13) at (cos{336},sin{336}) {}; 
    \draw (0) arc(0:308:1);
    \draw(0)--(9);\draw(2)--(4);\draw(3)--(12);\draw(8)--(10);
    \draw(5)--(7);\draw(6)--(13);
     
\end{tikzpicture}
\hspace{0mm}
\begin{tikzpicture}[scale=2.8]
\tikzstyle{vertex}=[circle, draw, fill=black, inner sep=0pt, minimum size=5pt]    
        \node[vertex, label=right:1](0) at (1,0) {};
    \node[vertex, label=above right:6](1) at (cos{26},sin{26}) {};
    \node[vertex, label=above right:6](2) at (cos{52},sin{52}) {};
    \node[vertex, label=above:3](3) at (cos{78},sin{78}) {};
    \node[vertex, label=above:5](4) at (cos{104}, sin{104}) {};
    \node[vertex, label=above left:5](5) at (cos{130}, sin{130}) {};
    \node[vertex, label=above left:2](6) at (cos{156}, sin{156}) {};
    \node[vertex, label=below left:4](7) at (cos{180}, sin{180}) {};
    \node[vertex, label=below left:4](8) at (cos{206}, sin{206}) {};
    \node[vertex, label=below:1](9) at (cos{232},sin{232}) {};
    \node[vertex, label=below:7](10) at (cos{258}, sin{258}) {};
    \node[vertex, label=below:7] (11) at (cos{284},sin{284}) {};    
    \node[vertex, label=below right:3] (12) at (cos{310},sin{310}) {}; 
    \node[vertex, label=below right:2] (13) at (cos{336},sin{336}) {};  
    \draw (0) arc(0:308:1);
    \draw(0)--(9);\draw(2)--(4);\draw(3)--(12);\draw(8)--(10);
    \draw(5)--(7);\draw(6)--(13);
     
\end{tikzpicture}

\caption {A square-free graph $\overline{KM_3}$ with a 6-clique-partition (left) and a frozen 7-clique-partition (right). Equivalently, a 6-colouring of the complement $KM_3$ of the graph shown (left) and a frozen 7-colouring of $KM_3$ (right).}
\label{fig:KM3}
\end{figure}
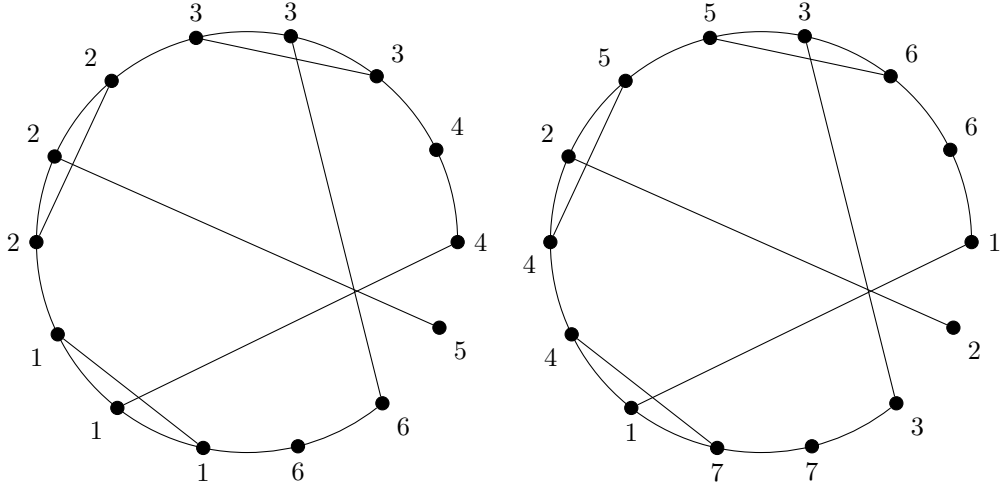


We now define a third class of graphs, ${KM}_q$ where $q \ge 2$. We obtain $\overline{{KM}_q}$ from $\overline{{ME}_q}$ by removing the edges of the path $u_1, u_2, \ldots, u_q$. Equivalently, we obtain ${KM}_q$ from ${ME}_q$ by adding the edges $u_1u_2, u_2u_3 \ldots, u_{q-1}u_q$ to  ${ME}_q$. See Figure \ref{fig:KM2} for $\overline{KM_2}$ and Figure \ref{fig:KM3} for $\overline{KM_3}$. 

\begin{theorem}\label{C4freeKMq} 
For $q \ge 2$, $\overline{{KM}_q}$ is $C_4$-free.
\end{theorem}

\begin{theorem}\label{thetaKMq}
For $q \ge 2$,
$\theta(\overline{{KM}_q})$ = $\alpha(\overline{{KM}_q})$ =
$2q$.
\end{theorem}

\begin{proof}
Let $q \ge 2$ be even. Create a clique partition of $\overline{{KM}_q}$ consisting of the following cliques:

\begin{itemize}
\item For $i=1, 2, \dots, q$, let the vertices of triangle $i$ be a clique in the clique partition
\item Let $\{u_0,u_1\}$ and $\{u_q,u_{q+1}\}$ be cliques of the clique partition
\item Let $\{u_2\},\{u_3\},\dots, \{u_{q-1}\}$ be cliques of the clique partition.
\end{itemize}

This clique partition has size $q+2+(q-2)=2q.$

Vertices $u_1, u_2, u_3, \dots, u_q,v_{11}, v_{21}, v_{31}, \ldots, v_{q1}$ form an independent set of size $2q$ in $\overline{{KM}_q}$. 
\end{proof}

\begin{theorem}\label{thm:frozencpKMq}
For $q \ge 2$, $\overline{{KM}_q}$ has a frozen $(2q+1)$-clique-partition.
\end{theorem}

\begin{proof}
The frozen clique partition of ${ME}_q$ given in Theorem \ref{thm:frozencp} is a frozen clique partition of ${KM}_q$.
\end{proof}

\begin{corollary} \label{KMcomplement}
For $q \ge 2$, ${KM}_q$ is a $2q$-chromatic $2K_2$-free graph with a frozen $(2q+1)$-colouring.
\end{corollary}

See Table \ref{tab:KMq graphs} for parameters of ${KM}_q$ graphs.

\begin{table} 
    \centering
    \begin{tabular}{|c|c|c|c|c|c|c|c|}
    \hline
        $q$ & $n$ & min  & max & \# edges & $\chi=\omega$ & \# colours in & (\# colours in  \\
         &  & degree & degree & &  & frozen colouring & frozen colouring) - $\chi$ 
         \\
         \hline
         & & & & & &  &     \\
         $q$ & $4q+2$ & $4q-2$ & $4q$ for $q > 2$ & $8q^2+q-2$ &  $2q$  & $2q+1$ &  $1$    \\
          &  & & $4q-1$ for $q=2$   &  &   &  &      \\
        & & & &  &  &  &   \\
         \hline
          & & & &  &  & &    \\
         2 & 10 & ~6 & ~7 & ~32  & ~4 & ~5 & 1  \\
         3 & 14 & 10 & 12 & ~63   & ~6  & ~7 & 1  \\
         4 & 18 & 14 & 16 & 130 & ~8  & ~9 & 1 \\
         5 & 22 & 18 & 20 & 203 & 10  & 11& 1  \\
         6 & 26 & 22 & 24 & 292 & 12  & 13 & 1 \\
         7 & 30 & 26 & 28 & 397 & 14  & 15 & 1  \\
         8 & 36 & 30 & 32 & 518 & 16  & 17 & 1  \\
         \hline    
    \end{tabular}
    \caption{Parameters of ${KM}_q$ graphs}
    \label{tab:KMq graphs}
\end{table}

We now define a fourth class of graphs. For  $q \ge 1$, $\overline{{KE}_q}$  is the graph with $6q$ vertices \\ 
$\{\{\cup \{v_{i1}, v_{i2}, v_{i3}\}: i=1,2,\ldots,2q\}$\\ whose edges are: 
\begin{itemize}
\item  the edges of a Hamiltonian cycle $C$: $v_{11}, v_{12}, v_{13}, v_{21}, v_{22}, v_{33},\ldots, v_{2q \mkern3mu 1}, v_{2q \mkern3mu 2}, v_{2q \mkern3mu 3}$
\item edges $v_{i2}v_{i+q \mkern3mu 2}$ for $i = 1, 2, \ldots, q$
\end{itemize}

As above, we refer to  $\{v_{i1}, v_{i2}, v_{i3}\}$ as \emph{triangle $i$}. Note that $\overline{{KE}_q}$ consists of a Hamiltonian cycle $C$ together with $2q$ edges which induce $2q$ vertex-disjoint triangles with consecutive pairs of edges of $C$, and $q$ more edges pairing the middle vertices $v_{i2}$ of ``opposite" triangles. The number of edges of $\overline{{KE}_q}$  is $9q$. 

\begin{figure}
\centering
\begin{tikzpicture}[scale=2.8]   
\tikzstyle{vertex}=[circle, draw, fill=black, inner sep=0pt, minimum size=5pt]    
        \draw (0,0) circle (1); 
         \node[vertex, label=right:4](0) at (1,0) {};
       \node[vertex, label=above:4](1) at (cos{20},sin{20}) {};
    \node[vertex, label=right:3](2) at (cos{40},sin{40}) {};
    \node[vertex, label=above:3](3) at (cos{60}, sin{60}) {};
    \node[vertex, label=above:3](4) at (cos{80}, sin{80}) {};
    \node[vertex, label=above:2](5) at (cos{100}, sin{100}) {};
    \node[vertex, label=above left:2](6) at (cos{120}, sin{120}) {};
    \node[vertex, label=above left:2](7) at (cos{140},sin{140}) {};
    \node[vertex, label=left:1](8) at (cos{160}, sin{160}) {};
    \node[vertex, label=left:1](9) at (cos{180},sin{180})
    {};
   \node[vertex, label=left:1](10) at (cos{200},sin{200}) {}; 
   \node[vertex, label=below left:6](11) at (cos{220},sin{220}) {}; 
   \node[vertex, label=below left:6](12) at (cos{240},sin{240}) {}; 
   \node[vertex, label=below:6](13) at (cos{260},sin{260}) {}; 
   \node[vertex, label=below:5](14) at (cos{280},sin{280}) {}; 
   \node[vertex, label=below:5](15) at (cos{300},sin{300}) {}; 
   \node[vertex, label=below right:5](16) at (cos{320},sin{320}) {}; 
   \node[vertex, label=below right:4](17) at (cos{340},sin{340}) {}; 
     \draw(0)--(9);\draw(1)--(17);\draw(2)--(4);\draw(3)--(12); \draw(5)--(7);\draw(6)--(15); \draw(8)--(10); \draw(11)--(13);
    \draw(14)--(16);
    
\end{tikzpicture}
\hspace{0mm}
\begin{tikzpicture}[scale=2.8]
\tikzstyle{vertex}=[circle, draw, fill=black, inner sep=0pt, minimum size=5pt]    
        \draw (0,0) circle (1); 
        \node[vertex, label=right:1](0) at (1,0) {};
    \node[vertex, label=right:6](1) at (cos{20},sin{20}) {};
    \node[vertex, label=above:6](2) at (cos{40},sin{40}) {};
    \node[vertex, label=above:3](3) at (cos{60}, sin{60}) {};
    \node[vertex, label=above:5](4) at (cos{80}, sin{80}) {};
    \node[vertex, label=above:5](5) at (cos{100}, sin{100}) {};
    \node[vertex, label=above left:2](6) at (cos{120}, sin{120}) {};
    \node[vertex, label=above left:4](7) at (cos{140},sin{140}) {};
    \node[vertex, label=left:4](8) at (cos{160}, sin{160}) {};
    \node[vertex, label=left:1](9) at (cos{180},sin{180})
    {};
   \node[vertex, label=left:9](10) at (cos{200},sin{200}) {}; 
   \node[vertex, label=below left:9](11) at (cos{220},sin{220}) {}; 
   \node[vertex, label=below left:3](12) at (cos{240},sin{240}) {}; 
   \node[vertex, label=below:8](13) at (cos{260},sin{260}) {}; 
   \node[vertex, label=below:8](14) at (cos{280},sin{280}) {}; 
   \node[vertex, label=below right:2](15) at (cos{300},sin{300}) {}; 
   \node[vertex, label=below right:7](16) at (cos{320},sin{320}) {}; 
   \node[vertex, label=below right:7](17) at (cos{340},sin{340}) {}; 
    \draw(0)--(9);\draw(1)--(17);\draw(2)--(4);\draw(3)--(12); \draw(5)--(7);\draw(6)--(15); \draw(8)--(10); \draw(11)--(13); \draw(14)--(16);
\end{tikzpicture}

\caption{A $C_4$-free graph $\overline{KE_3}$ with a 6-clique-partition (left) and a frozen 9-clique-partition (right). Equivalently, a 6-colouring of the complement $KE_3$ (left) and a frozen 9-colouring (right).}
\label{fig:KE3}
\end{figure}
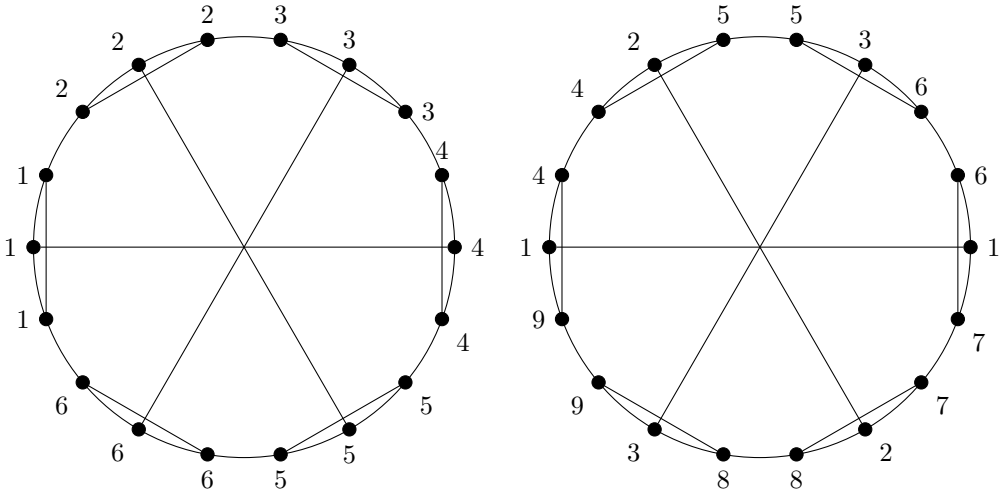

See Figure \ref{fig:KE2} for $\overline{KE_2}$  and Figure \ref{fig:KE3} for $\overline{KE_3}$. Note that $\overline{KE_1}$ is $\overline{C_6}$.

\begin{figure}
\centering
\begin{tikzpicture}[scale=2.8]   
\tikzstyle{vertex}=[circle, draw, fill=black, inner sep=0pt, minimum size=5pt]    
        \draw (0,0) circle (1); 
        \node[vertex, label=right:3](0) at (1,0) {};
    \node[vertex, label=above right:3](1) at (cos{30},sin{30}) {};
    \node[vertex, label=above:2](2) at (cos{60},sin{60}) {};
    \node[vertex, label=above:2](3) at (cos{90}, sin{90}) {};
    \node[vertex, label=above left:2](4) at (cos{120}, sin{120}) {};
    \node[vertex, label=left:1](5) at (cos{150}, sin{150}) {};
    \node[vertex, label=below left:1](6) at (cos{180}, sin{180}) {};
    \node[vertex, label=below:1](7) at (cos{210},sin{210}) {};
    \node[vertex, label=below:4](8) at (cos{240}, sin{240}) {};
    \node[vertex, label=below right:4](9) at (cos{270},sin{270})
    {};
   \node[vertex, label=below right:4](10) at (cos{300},sin{300}) {}; 
   \node[vertex, label=below right:3](11) at (cos{330},sin{330}) {}; 
    \draw(0)--(6);\draw(1)--(11);\draw(2)--(4);\draw(3)--(9); \draw(5)--(7);\draw(8)--(10);
    
\end{tikzpicture}
\hspace{0mm}
\begin{tikzpicture}[scale=2.8]
\tikzstyle{vertex}=[circle, draw, fill=black, inner sep=0pt, minimum size=5pt]    
        \draw (0,0) circle (1); 
        \node[vertex, label=right:1](0) at (1,0) {};
    \node[vertex, label=above right:4](1) at (cos{30},sin{30}) {};
    \node[vertex, label=above:4](2) at (cos{60},sin{60}) {};
    \node[vertex, label=above:2](3) at (cos{90}, sin{90}) {};
    \node[vertex, label=above left:3](4) at (cos{120}, sin{120}) {};
    \node[vertex, label=left:3](5) at (cos{150}, sin{150}) {};
    \node[vertex, label=below left:1](6) at (cos{180}, sin{180}) {};
    \node[vertex, label=below:6](7) at (cos{210},sin{210}) {};
    \node[vertex, label=below:6](8) at (cos{240}, sin{240}) {};
    \node[vertex, label=below right:2](9) at (cos{270},sin{270})
    {};
   \node[vertex, label=below right:5](10) at (cos{300},sin{300}) {}; 
   \node[vertex, label=below right:5](11) at (cos{330},sin{330}) {}; 
    \draw(0)--(6);\draw(1)--(11);\draw(2)--(4);\draw(3)--(9); \draw(5)--(7);\draw(8)--(10);
\end{tikzpicture}

\caption{A $C_4$-free graph $\overline{KE_2}$ with a 4-clique-partition (left) and a frozen 6-clique-partition (right). Equivalently, a 4-colouring of the complement $KE_2$  (left) and a frozen 6-colouring (right).}
\label{fig:KE2}
\end{figure}
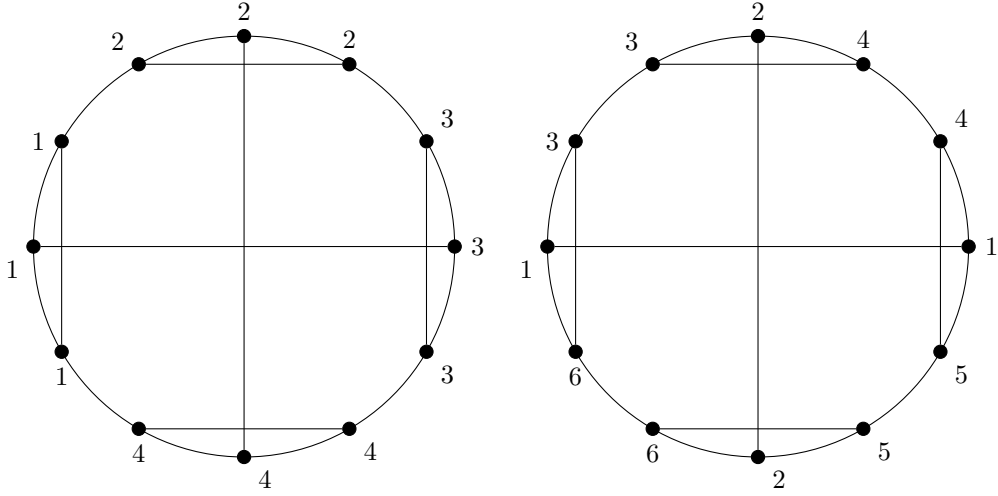

\begin{theorem}\label{KEC4free} 
For $q \ge 2$, $\overline{{KE}_q}$ is $C_4$-free.
\end{theorem}

\begin{theorem}\label{KEtheta}
For $q \ge 2$,
$\theta(\overline{{KE}_q}) = \alpha(\overline{{KE}_q}) = 2q$.
\end{theorem}

\begin{proof}
Let $q \ge 2$. Triangles $1, 2, \ldots, 2q$ form a clique partition of $\overline{{KE}_q}$. Vertices $v_{11}, v_{21}, \ldots, v_{2q \mkern3mu 1}$ form an independent set in $\overline{{ME}_q}$ graphs.
\end{proof}

\begin{theorem}\label{thm:KEfrozencp}
For $q \ge 2$, $\overline{{KE}_q}$ has a frozen $3q$-clique-partition.
\end{theorem}

\begin{proof}
Let $q \ge 2$. The following $3q$-clique-partition is frozen: \\$\mathcal{Q}= 
\{v_{12},v_{q+1 \mkern3mu 2} \}, 
\{v_{22},v_{q+2 \mkern3mu 2} \}, \ldots, 
\{v_{q2},v_{2q \mkern3mu 2} \},  
\{v_{13}, v_{21}\}, 
\{v_{23}, v_{31}\},\ldots,
\{v_{2q \mkern3mu 3}, v_{11}\}\}$.
\end{proof}

\begin{corollary} \label{KEcomplement}
For $q \ge 2$, ${KE}_q$ is a $2q$-chromatic $2K_2$-free graph with a frozen $3q$-colouring.
\end{corollary}

See Table \ref{tab:KEq graphs} for parameters of ${KE}_q$ graphs.

\begin{table} 
    \centering
    \begin{tabular}{|c|c|c|c|c|c|c|c|}
    \hline
        $q$ & $n$ & min  & max & \# edges & $\chi=\omega$ & \# colours in & (\# colours in  \\
         &  & degree & degree & &  & frozen colouring & frozen colouring) - $\chi$ 
         \\
         \hline
         & & & & & &  &     \\
         $q$ & $6q$ & $6q-4$ & $6q-4$ & $8q^2-12q$ &  $2q$  & $3q$ &  $q$    \\
        & & & &  &  &  &   \\
         \hline
          & & & &  &  & &    \\
         1 & 6  & ~2 & ~2 & ~6  & 2 & 3 & 1  \\
         2 & 12 & 8 & ~8 & 48   & 4  & 6 & 2  \\
         3 & 18 & 14 & 14 & 126 & 6  & 9 & 3 \\
         4 & 24 & 20 & 20 & 240 & 8  & 12& 4  \\
         5 & 30 & 26 & 26 & 390 & 10  & 15 & 5 \\
         6 & 36 & 32 & 32 & 576 & 12  & 18 & 6  \\
         7 & 42 & 38 & 38 & 798 & 14  & 21 & 7  \\
         8 & 48 & 44 & 44 & 1036& 16  & 24 & 8  \\
         \hline    
    \end{tabular}
    \caption{Parameters of ${KE}_q$ graphs}
    \label{tab:KEq graphs}
\end{table}

\begin{remark}
For any $q \ge 2$, one can obtain a $2q$-clique-partionable graph with a frozen $3q$-clique-partition by modifying the construction of $\overline{{KE}_q}$ as follows. Pair the vertices $\{v_{i2}: 1 \le i \le q\}$ in any way, and then join the members of each pair by an edge (rather than joining $v_{i2}$ to $v_{i+q \mkern3mu 2}$ as in the construction). To avoid creating a $C_4$, do not pair $v_{i2}$ with $v_{i+1  \mkern3mu 2}$ for $1 \le i \le q-1$ and do not pair $v_{iq}$ with $v_{i1}$.
\end{remark}

\section{An operation which preserves being \texorpdfstring{$2K_2$}\ -free and admitting a frozen colouring}
\label{sec:operation}

\begin{operation}
\label{op:subdivide}
Given a graph $H$ and adjacent vertices $x$ and $y$ in $H$, we subdivide the edge $xy$ to obtain a new graph $H^\prime$ by deleting the edge $xy$, adding two vertices $u$ and $v$, and adding edges $xu,~uv$, and $vy$; that is, the edge $xy$ is replaced by a path on four vertices: $x,u,v,y$.
\end{operation}

\begin{theorem}\label{thm:C_4 increase}
Let $H$ be a graph with a $k$-clique-partition $\mathcal{Q}$ and with a frozen $(k+1)$-clique-partition $\mathcal{F}$, and let $x$ and $y$ be adjacent vertices of $H$ which are in different cliques of $\mathcal{Q}$ such that either 

\begin{itemize}
\item [(1)] $x$ and $y$ are in different cliques of $\mathcal{F}$  or 
\item [(2)]	$\{x,y\}$ is a clique of $\mathcal{F}$.
\end{itemize}

\noindent Then the graph $H^\prime$ obtained by subdividing edge $xy$ as in Operation \ref{op:subdivide} is $(k+1)$-clique-partitionable and admits a frozen $(k+2)$-clique-partition.
\\Furthermore, 
\begin{itemize}
\item [(3)] if $\theta(H)=k$, then $\theta(H^\prime) = k+1$.
\item [(4)] if $H$ is $C_4$-free and if in case (1), $xy$ is not the middle edge of a diamond, then $H^\prime$ is $C_4$-free. 
\end{itemize}
\end{theorem}

\begin{proof}
Let $H$ be a graph with a $k$-clique-partition $\mathcal{Q}$ and with a frozen $(k+1)$-clique-partition $\mathcal{F}$, and let $x$ and $y$ be adjacent vertices of $H$ which are in different cliques of $\mathcal{Q}$. 
Let $H^\prime$ be the graph obtained by subdividing edge $xy$.\\

\begin{claim}
By adding $\{u,v\}$ to $\mathcal{Q}$ we obtain a $(k+1)$-clique-partition $\mathcal{Q^\prime}$ of $H^\prime$.    
\end{claim}

\begin{claim}
We can modify $\mathcal{F}$ to be a frozen colouring $\mathcal{F}^\prime$  of $H^\prime$ as follows. 

\begin{itemize}
    \item[]In Case (1): By adding $\{u,v\}$ to $\mathcal{F}$ we obtain a $(k+1)$-clique-partition $\mathcal{F^\prime}$ of $H^\prime$. 
    \item[]In Case (2): Remove $\{x,y\}$ from $\mathcal{F}$ and add $\{x,u\}$ and $\{v,y\}$ to obtain a $(k+1)$-clique-partition $\mathcal{F^\prime}$ of $H^\prime$.
\end{itemize}

\end{claim}

\begin{proof}

It is easy to see that $\mathcal{F}^\prime$ is a clique partition of $H^\prime$. We now prove that  $\mathcal{F}^\prime$ is frozen.

    In Case (1): In $H^\prime$, every vertex is nonadjacent to either $u$ or $v$ or both, so every vertex not in clique $\{u,v\}$ is nonadjacent to a vertex of $\{u,v\}$.

    Since $\mathcal{F}$ is a frozen clique partition of $H$, every vertex of $H$ is nonadjacent to some vertex of every clique of $\mathcal{F}$ other than the clique containing it, and this remains true when the edge $xy$ is deleted. 

    Thus, for every vertex $z$ of $H$ and every clique $Q$ of $\mathcal{F}^\prime$ other than the clique containing $z$, $z$ is nonadjacent to some vertex of $Q$.

    In any frozen clique partition, if there is a clique consisting of a single vertex, say $w$, then $w$ must be an isolated vertex. In $H$, $x$ and $y$ are adjacent, so neither is an isolated vertex, and thus there is vertex $x^\prime$ of $H$ different from $x$ in  the clique of $\mathcal{F}$ containing $x$ and a vertex $y^\prime$ different from $y$ in  the clique of $\mathcal{F}$ containing $y$. 

    Since $u$ is nonadjacent to every vertex of $H$ other than $x$, and in particular, is nonadjacent to $x^\prime$, it follows that $u$ is nonadjacent to some vertex of every clique of $\mathcal{F}^\prime$ other than $\{u,v\}$. Similarly, $v$ is nonadjacent to some vertex of every clique of $\mathcal{F}^\prime$ other than $\{u,v\}$.\\

    In Case (2): In $H^\prime$, vertex $u$ is nonadjacent to every vertex other than $x$ and $v$. Thus $u$ is nonadjacent to some vertex of every clique of $\mathcal{F}^\prime$ other than $\{x,u\}$. Analogously, $v$ is nonadjacent to some vertex of every clique of $\mathcal{F}^\prime$ other than $\{v,y\}$.

    Since $\mathcal{F}$ is a frozen clique partition of $H$, every vertex of $H$ is nonadjacent to some vertex of every clique of $\mathcal{F}$ other than the clique containing it. In particular, every vertex $w$ in $V(H)-\{x,y\}$ is nonadjacent to a vertex of each clique of $\mathcal{F} \setminus \{x,y\}$. Since $w$ is nonadjacent to $u$ and $v$, it follows $w$ is nonadjacent to some vertex of each clique of $\mathcal{F}^\prime = (\mathcal{F} \setminus \{x,y\}) \cup \{\{x,u\},\{v,y\}\}$. 
    
    Since $\mathcal{F}$ is a frozen clique partition of $H$, $x$ is nonadjacent to a vertex of every clique of $\mathcal{F}$ other than $\{x,y\}$.  Vertex $x$ is nonadjacent to $v \in \{v,y\} \in \mathcal{F}^\prime $. Thus vertex $x$ is nonadjacent to some vertex of every clique of  $\mathcal{F}^\prime$ other than $\{x,u\}$. Analogously, vertex $y$ is nonadjacent to some vertex of every clique of  $\mathcal{F}^\prime$ other than $\{v,y\}$. 
\end{proof}

\begin{claim}
If $\theta(H)=k$, then $\theta(H^\prime)=k+1$.
\end{claim}

\begin{proof} 
Assume $\theta(H)=k$. 

If there were a $(k-2)$-clique partition of $H-\{x,y\}$, then by adding $\{x,y\}$ to the clique partition, we would obtain a $(k-1)$-clique-partition of $H$, which is a contradiction. So $\theta(H-\{x,y\}) \ge k-1$.

By Claim 1, $\theta(H^\prime) \le k+1$. We need to show that there is no $k$-clique-partition of $H^\prime$. First, consider a clique partition of $H^\prime$ where $u$ and $v$ are in different cliques. Since $u$ and $v$ are each anticomplete to 
$H-\{x,y\}$ and $\theta(H-\{x,y\}) \ge k-1$, a total of at least $k+1$ cliques would be required. Now consider a clique partition of $H^\prime$ where $u$ and $v$ are in the same clique. This clique must then be $\{u,v\}$, and thus the clique partition must have at least $\theta(H)+1 = k+1$ cliques.
\end{proof}

\begin{claim}
If $H$ is $C_4$-free, then $H^\prime$ is $C_4$-free.
\end{claim}

\begin{proof} 
Assume $H$ is $C_4$-free.

In $H^\prime$, $u$ and $v$ are adjacent and each have degree 2, so any $C_4$ containing one of them, must contain the other, and then also contain $u$'s only other neighbour, which is $x$, and $v$'s only other neighbour, which is $y$, but $xy$ is not an edge of $H^\prime$, so no such $C_4$ exists.

In constructing $H^\prime$ from $H$, the edge $xy$ is removed. This could create a $C_4$ if $xy$ was the middle edge of a diamond in $H$. This is excluded by hypothesis in Case (1). In Case (2), $\{x,y\}$ is a clique in the frozen clique partition $\mathcal{F}$. If there were a vertex $w$ adjacent to both $x$ and $y$ in $H$, then $\mathcal{F}$ would not be frozen. Thus $xy$ cannot be the middle edge of a diamond in $H$.
\end{proof}

\end{proof}

\section{k-chromatic \texorpdfstring{$2K_2$}\ -free graphs which admit a frozen (k+1)-colouring for all k \texorpdfstring{$\ge$}\ 4}
\label{sec:all}

\begin{figure}
\centering
\begin{tikzpicture}[scale=2.8]   
\tikzstyle{vertex}=[circle, draw, fill=black, inner sep=0pt, minimum size=5pt]    
        \draw (0,0) circle (1); 
        \node[vertex, label=right:2](0) at (1,0) {};
    \node[vertex, label=above right:2](1) at (cos{36},sin{36}) {};
    \node[vertex, label=above:2](2) at (cos{72},sin{72}) {};
    \node[vertex, label=above:1](3) at (cos{108}, sin{108}) {};
    \node[vertex, label=above left:1](4) at (cos{144}, sin{144}) {};
    \node[vertex, label=left:1](5) at (cos{180}, sin{180}) {};
    \node[vertex, label=below left:5](6) at (cos{206}, sin{206}) {};
    \node[vertex, label=below:5](7) at (cos{232},sin{232}) {};
    \node[vertex, label=below:4](8) at (cos{258}, sin{258}) {};
    \node[vertex, label=below right:4](9) at (cos{284},sin{284})
    {};
   \node[vertex, label=below right:3](10) at (cos{310},sin{310}) {}; 
   \node[vertex, label=below right:3](11) at (cos{336},sin{336}) {}; 
    \draw(0)--(2);\draw(1)--(7);\draw(3)--(5);\draw(4)--(10);
    
\end{tikzpicture}
\hspace{0mm}
\begin{tikzpicture}[scale=2.8]
\tikzstyle{vertex}=[circle, draw, fill=black, inner sep=0pt, minimum size=5pt]    
        \draw (0,0) circle (1); 
        \node[vertex, label=right:4](0) at (1,0) {};
    \node[vertex, label=above right:2](1) at (cos{36},sin{36}) {};
    \node[vertex, label=above:3](2) at (cos{72},sin{72}) {};
    \node[vertex, label=above:3](3) at (cos{108}, sin{108}) {};
    \node[vertex, label=above left:1](4) at (cos{144}, sin{144}) {};
    \node[vertex, label=left:6](5) at (cos{180}, sin{180}) {};
    \node[vertex, label=below left:6](6) at (cos{206}, sin{206}) {};
    \node[vertex, label=below:2](7) at (cos{232},sin{232}) {};
    \node[vertex, label=below:5](8) at (cos{258}, sin{258}) {};
    \node[vertex, label=below right:5](9) at (cos{284},sin{284})
    {};
   \node[vertex, label=below right:1](10) at (cos{310},sin{310}) {}; 
   \node[vertex, label=below right:4](11) at (cos{336},sin{336}) {}; 
    \draw(0)--(2);\draw(1)--(7);\draw(3)--(5);\draw(4)--(10);
\end{tikzpicture}

\caption{A $C_4$-free graph with a 5-clique-partition (left) and a frozen 6-clique-partition (right). Equivalently, a 5-colouring of the complement (left) and a frozen 6-colouring (right).}
\label{fig:5c6f}

\begin{tikzpicture}[scale=2.8]   
\centering
\tikzstyle{vertex}=[circle, draw, fill=black, inner sep=0pt, minimum size=5pt]    
        \draw (0,0) circle (1); 
        \node[vertex, label=right:2](0) at (1,0) {};
    \node[vertex, label=above right:2](1) at (cos{36},sin{36}) {};
    \node[vertex, label=above:2](2) at (cos{72},sin{72}) {};
    \node[vertex, label=above:1](3) at (cos{108}, sin{108}) {};
    \node[vertex, label=above left:1](4) at (cos{144}, sin{144}) {};
    \node[vertex, label=left:1](5) at (cos{180}, sin{180}) {};
    \node[vertex, label=below left:7](6) at (cos{197}, sin{197}) {};
    \node[vertex, label=below left:7](7) at (cos{214},sin{214}) {};
    \node[vertex, label=below left:6](8) at (cos{230}, sin{230}) {};
    \node[vertex, label=below left:6](9) at (cos{246},sin{246})
    {};
   \node[vertex, label=below:5](10) at (cos{262},sin{262}) {}; 
   \node[vertex, label=below:5](11) at (cos{278},sin{278}) {}; 
    \node[vertex, label=below right:4](12) at (cos{294},sin{294}) {}; 
     \node[vertex, label=below right:4](13) at (cos{310},sin{310}) {}; 
      \node[vertex, label=below right:3](14) at (cos{326},sin{326}) {}; 
       \node[vertex, label=below right:3](15) at (cos{343},sin{343}) {}; 
    \draw(0)--(2);\draw(1)--(7);\draw(3)--(5);\draw(4)--(14);
\end{tikzpicture}
\hspace{0mm}
\begin{tikzpicture}[scale=2.8]
\tikzstyle{vertex}=[circle, draw, fill=black, inner sep=0pt, minimum size=5pt]      
\draw (0,0) circle (1); 
        \node[vertex, label=right:4](0) at (1,0) {};
    \node[vertex, label=above right:2](1) at (cos{36},sin{36}) {};
    \node[vertex, label=above:3](2) at (cos{72},sin{72}) {};
    \node[vertex, label=above:3](3) at (cos{108}, sin{108}) {};
    \node[vertex, label=above left:1](4) at (cos{144}, sin{144}) {};
    \node[vertex, label=left:8](5) at (cos{180}, sin{180}) {};
    \node[vertex, label=below left:8](6) at (cos{197}, sin{197}) {};
    \node[vertex, label=below left:2](7) at (cos{214},sin{214}) {};
    \node[vertex, label=below left:7](8) at (cos{230}, sin{230}) {};
    \node[vertex, label=below left:7](9) at (cos{246},sin{246})
    {};
   \node[vertex, label=below:6](10) at (cos{262},sin{262}) {}; 
   \node[vertex, label=below:6](11) at (cos{278},sin{278}) {}; 
    \node[vertex, label=below right:5](12) at (cos{294},sin{294}) {}; 
     \node[vertex, label=below right:5](13) at (cos{310},sin{310}) {}; 
      \node[vertex, label=below right:1](14) at (cos{326},sin{326}) {}; 
       \node[vertex, label=below right:4](15) at (cos{343},sin{343}) {}; 
    \draw(0)--(2);\draw(1)--(7);\draw(3)--(5);\draw(4)--(14);
\end{tikzpicture}

\caption{A $C_4$-free graph with a 7-clique-partition (left) and a frozen 8-clique-partition (right). Equivalently, a 7-colouring of the complement (left) and a frozen 8-colouring (right).}
\label{fig:7c8f}
\end{figure}   

\begin{theorem}\label{C4 all k}
For every $k \ge 4$, there is a $C_4$-free graph with clique partition number $k$ which admits a frozen $(k+1)$-clique partition.
\end{theorem}

\begin{proof}
One way to construct the graphs described in the theorem is to start with $\overline{ME_2}$ which is a $C_4$-free graph with clique partition number $4$ and with a frozen 5-clique-partition and then apply the operation described in Theorem \ref{thm:C_4 increase} with $x=u_1$ and $y=u_2$. These two vertices are in different cliques in both the 4-clique-partition and in the frozen 5-clique-partition, so Case (1) will be applied. The additional hypothesis holds in this case. The result is a $C_4$-free graph with clique partition number $5$ and with a frozen 6-clique-partition. Note that the two added vertices are a clique of size 2 in both the 5-clique-partition and the frozen 6-clique-partition. One can then apply the operation again, with $x=u_1$ and $y$ being the vertex $u$ of the previous operation to obtain a $C_4$-free graph with clique partition number $6$ and with a frozen 7-clique-partition. One can continue this process, always choosing $x=u_1$ and $y$ being the vertex $u$ of the previous operation. This class of graphs is illustrated in Figures \ref{fig:5c6f} and \ref{fig:7c8f} and can be described as follows: For $t \ge 4$, to obtain a $C_4$-free graph with clique partition number $t$ and with a frozen $(t+1)$-clique-partition, start with $\overline{ME_2}$ and subdivide the edge $u_1u_2$ by $2(t-4)$ vertices (in other words, replace the edge $u_1u_2$ by a path $u_1, w_1, w_2, \ldots, w_{2t-9}, w_{2t-8}, u_2$).    
\end{proof}

\begin{corollary}
\label{cor:2K2 all k}
For every $k \ge 4$, there is a $k$-chromatic $2K_2$-free graph with a frozen $(k+1)$-colouring.
\end{corollary}

\begin{remark}
 There are many other ways to apply the operation described in the proof of Theorem \ref{thm:C_4 increase} to prove Theorem \ref{C4 all k} - it is not necessary to choose the same vertices as $x$ and $y$ as above. Since in the frozen clique partitions of $\overline{{ME}_q}$ and of $\overline{{{ME}_q}^*}$ given in Theorems \ref{thm:frozencp} and \ref{thm:frozen-}, all cliques have size 2, either of the operations from the proof of Theorem \ref{thm:C_4 increase} can be used. 
\end{remark}

Here is the operation described directly for colourings.

\begin{operation}
\label{op:2K2}
Given a graph $G$ and nonadjacent vertices $x$ and $y$ in $G$, we define the following operation to create a new graph $G^\prime$. Define $G^\prime$ to be the graph $G$ together with two additional vertices  $u$ and $v$ and with edges $vx,~ xy$ and $yu$; join $u$ and $v$ to all vertices of $G-\{x,y\}$.
\end{operation}

\begin{corollary}\label{cor:2K_2 increase}
Let $G$ be a $k$-colourable graph with a $k$-colouring $\beta$ and a frozen $(k+1)$-colouring $\gamma$, and let $x$ and $y$ be nonadjacent vertices of $G$ such that $\beta(x) \neq \beta(y)$ and such that either 

\begin{itemize}
\item [(1)] $\gamma(x) \neq \gamma(y)$, or 
\item [(2)]	$\{x,y\}$ is a colour class of $\gamma$.
\end{itemize}

\noindent Then the graph $G^\prime$ of Operation \ref{op:2K2} is $(k+1)$-colourable and admits a frozen $(k+2)$-colouring.
\\Furthermore, 
\begin{itemize}
\item [(3)] if $G$ is $k$-chromatic, then $G^\prime$ is $(k+1)$-chromatic.
\item [(4)] if $G$ is $2K_2$-free and if in case (1), there is no edge $rs$ such that 
$\{r,s\}$ is anticomplete to $\{x,y\}$, then $G^\prime$ is $2K_2$-free. 
\end{itemize}
\end{corollary}

\section{Some curiosities and open problems}
\label{sec:curious}

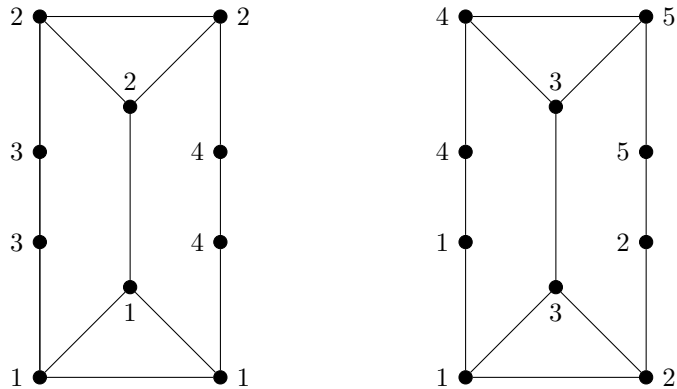
\begin{figure}[ht]
\centering
\begin{tikzpicture}[scale=1.2]
\tikzstyle{vertex}=[circle, draw, fill=black, inner sep=0pt, minimum size=5pt]
    \node[vertex, label=left:1](0) at (0,0) {};
    \node[vertex, label=right:1](1) at (2,0) {};
    \node[vertex, label=below:1](2) at (1,1) {};
    \node[vertex, label=left:2](3) at (0,4) {};
    \node[vertex, label=above:2](4) at (1,3) {};
    \node[vertex, label=right:2](5) at (2,4) {};
    \node[vertex, label=left:3](6) at (0,1.5) {};
    \node[vertex, label=left:3](7) at (0,2.5) {};
    \node[vertex, label=left:4](8) at (2,1.5) {};
    \node[vertex, label=left:4](9) at (2,2.5) {};
    \draw(0)--(2);\draw(0)--(1);\draw(0)--(3);\draw(2)--(1);\draw(1)--(8);\draw(8)--(9);\draw(0)--(6);\draw(6)--(7);\draw(3)--(7);\draw(5)--(3);\draw(3)--(4);\draw(4)--(5);\draw(9)--(5);\draw(2)--(4);
\end{tikzpicture}
\hspace{20mm}
\begin{tikzpicture}[scale=1.2]
\tikzstyle{vertex}=[circle, draw, fill=black, inner sep=0pt, minimum size=5pt]    
    \node[vertex, label=left:1](0) at (0,0) {};
    \node[vertex, label=right:2](1) at (2,0) {};
    \node[vertex, label=below:3](2) at (1,1) {};
    \node[vertex, label=left:4](3) at (0,4) {};
    \node[vertex, label=above:3](4) at (1,3) {};
    \node[vertex, label=right:5](5) at (2,4) {};
    \node[vertex, label=left:1](6) at (0,1.5) {};
    \node[vertex, label=left:4](7) at (0,2.5) {};
    \node[vertex, label=left:2](8) at (2,1.5) {};
    \node[vertex, label=left:5](9) at (2,2.5) {};
    \draw(0)--(2);\draw(0)--(1);\draw(0)--(6);\draw(2)--(1);\draw(1)--(8);\draw(8)--(9);\draw(2)--(4);\draw(3)--(7);\draw(6)--(7);\draw(5)--(3);\draw(3)--(4);\draw(4)--(5);\draw(9)--(5);
\end{tikzpicture}
\caption{A 4-clique-partition of 
$\overline{KM_2} \cong \overline{H_3}$
(left) and a frozen 5-clique-partition (right).}
\label{fig:coC6''}
\end{figure}

Subdividing an edge of a $C_4$ gets rid of that $C_4$. The complement $\overline{C_6}$ of $C_6$ contains three $C_4$s; each pair of $C_4$s intersect in a distinct edge. By applying Operation \ref{op:subdivide} to  two of these three edges, we obtain 
$\overline{KM_2}$ 
which is $C_4$-free. See Figure \ref{fig:coC6''}. Thus besides preserving $2K_2$-freeness of a graph, our operation can transform a graph containing $2K_2$s into a $2K_2$-free graph.

\begin{figure}[ht]
\centering
\begin{tikzpicture}[scale=1.2]
\tikzstyle{vertex}=[circle, draw, fill=black, inner sep=0pt, minimum size=5pt]
    \node[vertex, label=left:1](0) at (0,0) {};
    \node[vertex, label=right:1](1) at (3,0) {};
    \node[vertex, label=below:1](2) at (1,1) {};
    \node[vertex, label=below:1](3) at (2,1) {};
    \node[vertex, label=left:3](4) at (0,2) {};
    \node[vertex, label=left:4](5) at (1,2) {};
    \node[vertex, label=right:5](6) at (2,2) {};
    \node[vertex, label=left:3](7) at (0,2.5) {};
    \node[vertex, label=left:4](8) at (1,2.5) {};
    \node[vertex, label=right:5](9) at (2,2.5) {};
    \node[vertex, label=above:2](10) at (1,3.5) {};
    \node[vertex, label=above:2](11) at (2,3.5) {};
    \node[vertex, label=left:2](12) at (0,4.5) {};
    \node[vertex, label=right:2](13) at (3,4.5) {};
    
    \draw(0)--(1);\draw(0)--(2);\draw(0)--(3); \draw(0)--(4);\draw(1)--(2);\draw(1)--(3);\draw(1)--(13);\draw(2)--(3);\draw(2)--(5);\draw(3)--(6);\draw(4)--(7);\draw(5)--(8);\draw(6)--(9);\draw(7)--(12);\draw(8)--(10);\draw(9)--(11);\draw(10)--(11); \draw(10)--(12);\draw(10)--(13);\draw(11)--(12);\draw(11)--(13);\draw(12)--(13);
\end{tikzpicture}
\hspace{20mm}
\begin{tikzpicture}[scale=1.2]
\tikzstyle{vertex}=[circle, draw, fill=black, inner sep=0pt, minimum size=5pt]    
    \node[vertex, label=left:1](0) at (0,0) {};
    \node[vertex, label=right:4](1) at (3,0) {};
    \node[vertex, label=below:2](2) at (1,1) {};
    \node[vertex, label=below:3](3) at (2,1) {};
    \node[vertex, label=left:1](4) at (0,2) {};
    \node[vertex, label=left:2](5) at (1,2) {};
    \node[vertex, label=right:3](6) at (2,2) {};
    \node[vertex, label=left:5](7) at (0,2.5) {};
    \node[vertex, label=left:6](8) at (1,2.5) {};
    \node[vertex, label=right:7](9) at (2,2.5) {};
    \node[vertex, label=above:6](10) at (1,3.5) {};
    \node[vertex, label=above:7](11) at (2,3.5) {};
    \node[vertex, label=left:5](12) at (0,4.5) {};
    \node[vertex, label=right:4](13) at (3,4.5) {};
    
    \draw(0)--(1);\draw(0)--(2);\draw(0)--(3); \draw(0)--(4);\draw(1)--(2);\draw(1)--(3);\draw(1)--(13);\draw(2)--(3);\draw(2)--(5);\draw(3)--(6);\draw(4)--(7);\draw(5)--(8);\draw(6)--(9);\draw(7)--(12);\draw(8)--(10);\draw(9)--(11);\draw(10)--(11); \draw(10)--(12);\draw(10)--(13);\draw(11)--(12);\draw(11)--(13);\draw(12)--(13);
\end{tikzpicture}
\caption{A 5-clique-partition of  
$\overline{H_4}$
(left) and a frozen 7-clique-partition (right).}
\label{fig:H4complement}
\end{figure}
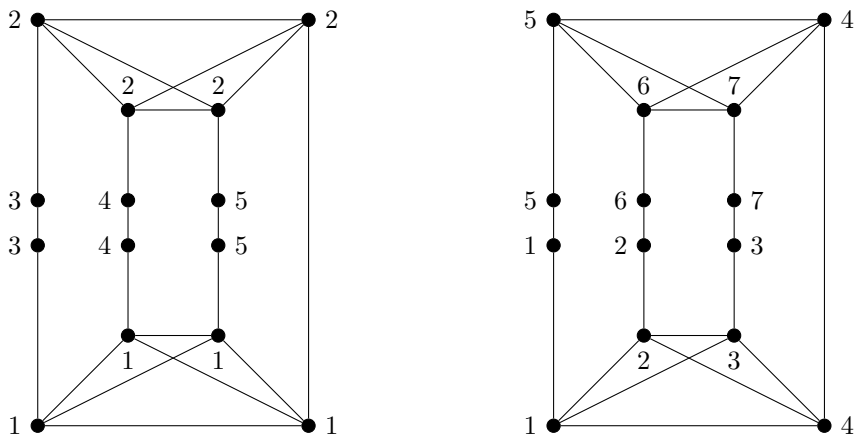

Recall that for $t \ge 2$, the graph $B_t$ is  $K_{t,t}$ with a perfect matching removed, and is 2-chromatic and admits a frozen $t$-colouring. The complement of $B_t$ consists of two copies of $K_t$ with a perfect matching $M_t$ joining each vertex of one copy to a distinct vertex of the other copy. Note that $B_t$ contains many $2K_2$s and (equivalently) $\overline{B_t}$ contains many $C_4$s.  By applying Operation \ref{op:subdivide} to all but one edge of $M_t$ in $\overline{B_t}$ where $t \ge 3$, we obtain a $C_4$-free graph $\overline{H_t}$ which is $(t+1)$-clique-partitionable and admits a frozen $(2t-1)$-clique partition. Note that $\overline{H_3}$ is isomorphic to $\overline{KM_2}.$ See Figure \ref{fig:coC6''} for $\overline{H_3}$ and Figure \ref{fig:H4complement} for $\overline{H_4}$.

\begin{figure}[ht]
\centering
\begin{tikzpicture}[scale=2.7]
\tikzstyle{vertex}=[circle, draw, fill=black, inner sep=0pt, minimum size=5pt]
    
    \node[vertex, label=above:3](0) at (0,1) {};
    \node[vertex, label=above:3](1) at (cos{67.5},sin{67.5}) {};
    \node[vertex, label=above:3](2) at (cos{45},sin{45}) {};
    \node[vertex, label=above right:7](3) at (cos{22.5}, sin{22.5}) {};
    \node[vertex, label=right:7](4) at (1,0) {};
    \node[vertex, label=below right:6](5) at (cos{22.5}, -sin{22.5}) {};
    \node[vertex, label=below right:1](6) at (cos{45},-sin{45}) {};
    \node[vertex, label=below:1](7) at (cos{67.5},-sin{67.5}) {};
    \node[vertex, label=below:4](8) at (0,-1) {};
    \node[vertex, label=below left:4](9) at (-cos{67.5},-sin{67.5}) {};
     \node[vertex, label=below left:4](10) at (-cos{45},-sin{45}) {};
    \node[vertex, label=left:2](11) at (-cos{22.5}, -sin{22.5}) {};
    \node[vertex, label=left:2](12) at (-1,0) {};
    \node[vertex, label=above left:2](13) at (-cos{22.5},sin{22.5}) {};
    \node[vertex, label=above left:5](14) at (-cos{45},sin{45}) {};
    \node[vertex, label=above:5](15) at (-cos{67.5},sin{67.5}) {};
    \draw (4) arc(0:315:1);
\draw(0)--(2); \draw(1)--(6);\draw(4)--(9);\draw(8)--(10);\draw(11)--(13);\draw(5)--(12);
    
\end{tikzpicture}
\hspace{0mm}
\begin{tikzpicture}[scale=2.7]
\tikzstyle{vertex}=[circle, draw, fill=black, inner sep=0pt, minimum size=5pt]    
    \node[vertex, label=above:5](0) at (0,1) {};
    \node[vertex, label=above:3](1) at (cos{67.5},sin{67.5}) {};
    \node[vertex, label=above:8](2) at (cos{45},sin{45}) {};
    \node[vertex, label=above right:8](3) at (cos{22.5}, sin{22.5}) {};
    \node[vertex, label=right:6](4) at (1,0) {};
    \node[vertex, label=below right:7](5) at (cos{22.5}, -sin{22.5}) {};
    \node[vertex, label=below right:3](6) at (cos{45},-sin{45}) {};
    \node[vertex, label=below:1](7) at (cos{67.5},-sin{67.5}) {};
    \node[vertex, label=below:1](8) at (0,-1) {};
    \node[vertex, label= below left:6](9) at (-cos{67.5},-sin{67.5}) {};
     \node[vertex, label=left:4](10) at (-cos{45},-sin{45}) {};
    \node[vertex, label=left:4](11) at (-cos{22.5}, -sin{22.5}) {};
    \node[vertex, label=left:7](12) at (-1,0) {};
    \node[vertex, label=above left:2](13) at (-cos{22.5},sin{22.5}) {};
    \node[vertex, label=above:2](14) at (-cos{45},sin{45}) {};
    \node[vertex, label=above:5](15) at (-cos{67.5},sin{67.5}) {};
    
    \draw (4) arc(0:315:1);
\draw(0)--(2); \draw(1)--(6);\draw(4)--(9);\draw(8)--(10);\draw(11)--(13);\draw(5)--(12);
\end{tikzpicture}
\caption{A 7-clique-partition (left) and a frozen 8-clique-partition (right) of a $C_4$-free graph. Equivalently, a 7-colouring (left) and a frozen 8-colouring (right) of the complement \cite{feghali2021}.}
\label{fig:OwenCarlComplement}
\end{figure}
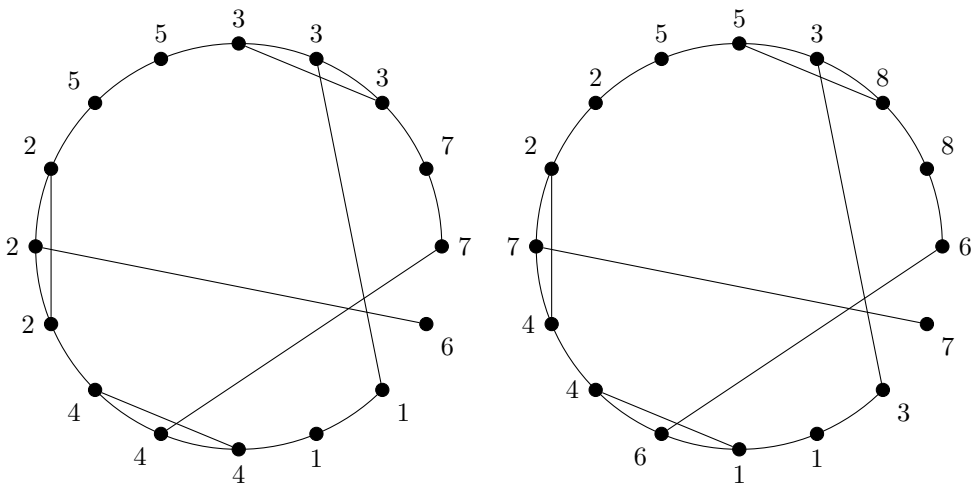

In Figure \ref{fig:OwenCarlComplement} is the complement of the $2K_2$-free graph given by Feghali and Merkel in \cite{feghali2021} with their 7-colouring (shown as a 7-clique-partition) and their frozen 8-colouring (shown as a frozen 8-clique-partition). The complement of their graph is very similar to our $\overline{KM_3}$. In fact, the complement of their graph is $\overline{KM_3}$ with Operation \ref{op:subdivide} applied once (to get the vertices in clique 5 of the 7-clique-partition).

We conclude with two open problems:
\begin{itemize}
    \item The remaining case of Question \ref{question}: Does there exist a $3$-colourable $P_5$-free graph with a frozen $4$-colouring? 
    \item The remaining case for a dichotomy theorem for recolouring graphs where two 4-vertex graphs are forbidden as induced subgraphs: Is the class of 4-chromatic $(2K_2, K_4)$-free graphs which contain a triangle recolourable?
\end{itemize}



\end{document}